\documentclass[preprint]{elsarticle}
\usepackage{soul,color}
\usepackage{amsmath,amsfonts,amssymb}
\usepackage[utf8]{inputenc}
\usepackage[T1]{fontenc}
\usepackage{units}
\usepackage{hyperref}
\usepackage{float}
\usepackage{subcaption}
\usepackage{fullpage}

\journal{Engineering Analysis with Boundary Elements}

\newcommand{\dpar}[2]{\ensuremath{\frac{\partial #1}{\partial #2}}}
\newcommand{\T}{\mathsf{T}}
\renewcommand{\b}{\boldsymbol}
\newcommand{\vp}{\vec{p}}
\renewcommand{\L}{\mathcal{L}}

\makeatletter
\def\ps@pprintTitle{%
  \let\@oddhead\@empty
  \let\@evenhead\@empty
  \def\@oddfoot{
Received 4 August 2017; Received in revised form 4 November 2017; Accepted 1 
January 2018. \hfill }%
  \let\@evenfoot\@oddfoot}
\makeatother

\begin{document}

\begin{frontmatter}

\title{Refined Meshless Local Strong Form solution of \\ Cauchy-Navier equation on an irregular domain}

\author{J. Slak}
\author{G. Kosec$^*$}
\address{``Jožef Stefan'' Institute, Parallel and Distributed systems 
Laboratory, Jamova 39, 1000 Ljubljana, Slovenia}

\cortext[G. Kosec]{Corresponding author}
\ead{gregor.kosec@ijs.si}

\begin{abstract}

This paper considers a numerical solution of a linear elasticity problem, namely
the Cauchy-Navier equation, using a strong form method based on a local Weighted Least Squares (WLS)
approximation. The main advantage of the employed numerical approach, also
referred to as a Meshless Local Strong Form method, is its generality in terms of
approximation setup and positions of computational nodes. In this paper,
flexibility regarding the nodal position is demonstrated through two numerical
examples, i.e.~a drilled cantilever beam, where an irregular domain is treated
with a relatively simple nodal positioning algorithm, and a Hertzian contact
problem, where again, a relatively simple h-refinement algorithm is used to
extensively refine discretization under the contact area. The results are
presented in terms of accuracy and convergence rates, using
different approximations and refinement setups, namely Gaussian and
monomial based approximations, and a comparison of execution time for each block of
the solution procedure.

\end{abstract}


\end{frontmatter}


\section{Introduction}
Linear elasticity problems, governed by the Cauchy-Navier equation, are typically
addressed in their weak form with the Finite Elements Method (FEM)~\cite{fem_solid}.
However, the problem has also been addressed in its strong form in the past,
e.g.\ component-wise iterative solution with the Finite Differences Method
(FDM)~\cite{volume_based_finite_difference}
and with the Finite Volumes Method (FVM)~\cite{solid_FVM1}.
Besides mesh based methods, meshless methods have also been employed for solving
solid mechanics problems in strong and weak form~\cite{meshless_solid,
strong_form_solid}.
The conceptual difference between meshless methods and mesh based methods is in
the treatment of relations between nodes. In the mesh based methods the nodes
need to be structured into polygons (mesh) that cover the whole computational
domain, while on the other hand, meshless methods fully define relations between
nodes through the relative inter nodal positions~\cite{trobec_kosec}, with an
immediate consequence of greater generality of the meshless methods.

Strong form meshless methods can be understood as generalizations of FDM, where
instead of predetermined interpolation over a local support, a more general
approach with variable support and basis is used to evaluate partial
differential operators~\cite{Kouszeg1}, e.g.\ collocation using
Radial Basis Functions~\cite{strong_form_solid,Kouszeg_PN}
or approximation with monomial basis ~\cite{Prax}. There are many
other methods with more or less similar methodology introducing new variants of
the strong form meshless principle~\cite{meshless_review}.  On the other hand,
weak form meshless methods are generalizations of FEM. Probably the most known
method among weak form meshless methods is the Meshless Local Petrov Galerkin
Method (MLPG)~\cite{Atluri}, where for each integration point a local support is
used to evaluate field values and weight functions of a Moving Least Squares
(MLS) approximation are used as test functions. In last few decades there
have been many variants of MLPG introduced to mitigate numerical instabilities
and to improve accuracy and convergence rate, etc.~\cite{meshless_review}.

In general, recent developments in meshless community are vivid, ranging from
analyses of computer execution on different platforms~\cite{trobec_kosec,
meshless_gpu}, reducing computational cost by introducing a piecewise
approximation~\cite{li_piece-wise_2017} to implementation of more complex
multi-phase flow~\cite{kovarik_meshless_2016}, and many more.

This paper extends the spectra of published papers with a generalized
formulation of a local strong form meshless method, termed Meshless Local Strong
Form Method (MLSM) enriched with h-refinement~\cite{Kouszeg_cmc} and ability to
discretize arbitrary domains~\cite{Kouszeg1}.

The introduced meshless approach is demonstrated on a solution of a benchmark
cantilever beam case~\cite{slaughter2012linearized} and a Hertzian contact
problem~\cite{pereira2016convergence}. The results are presented in terms of
displacement and stress plots, comparison against closed form solutions,
convergence analyses, and execution time analyses.

The goal of this paper is to demonstrate generality of MLSM that is driven by
the fact that all the building blocks of the method depend only on the relative
positions between the computational nodes. This is a very useful feature,
especially when dealing with problems in multidimensional spaces, complex
geometries, and moving boundaries. This feature can be also exploited to write
elegant generic code~\cite{utils_web}.

The rest of the paper is organized as follows: in section~\ref{sec:mlsm} the
MLSM principle is explained, in section~\ref{sec:governing-problem} the
governing problem is introduced, section~\ref{sec:solution-procedure} is focused
on solution procedure, section~\ref{sec:results} focuses on discussing
the results, and finally, the paper offers some conclusions and guidelines for
future work in the last section.

\section{MLSM formulation}
\label{sec:mlsm}
The core of the spatial discretization used in this paper is a local
approximation of a considered field over the overlapping local support domains,
i.e.\ in each node we use approximation over a small local subset of neighbouring
$n$ nodes. The trial function $\hat u$ is thus introduced as
\begin{equation}
  \hat{u}(\vp)=\sum_{i=1}^m \alpha_i b_i(\vp)=\b{b}(\vp)^\T \b\alpha,
  \label{eq:u-definition}
\end{equation}
with $m$, $\b{\alpha }$, $\b{b}$ and $\vec{p}$ standing for the number of basis
functions, approximation coefficients, basis functions and the position vector,
respectively.  In cases when the number of basis functions and the number of
nodes in the support domain are the same, $n=m$, the determination of
coefficients $\b{\alpha }$ simplifies to solving a system of $n$ linear
equations, resulting from evaluating equation~\eqref{eq:u-definition} in each
support node and setting it equal to a true value $u(\vp_j)$, for $j$ from 1 to $n$:
\begin{equation}
  u_j := u(\vp_j) = \b{b}(\vp_j)^\T \b{\alpha},
\end{equation}
where $\vp_j$ are positions of support nodes and $u_j$ is the actual value of
considered field in the support node $\vp_j$. The above system can be written in
matrix form as
\begin{equation}
  \b{u}=\b{B}\b{\alpha},
  \label{eq:basic-system}
\end{equation}
where $\b{B}$ stands for coefficient matrix with elements $B_{ji}=b_i(\vp_j)$.
The most known method that uses such an approach is the Local Radial Basis
Function Collocation Method (LRBFCM) that has been recently
used in various problems~\cite{strong_form_solid,Kouszeg_PN}.

In cases when the number of support nodes is higher than the number of basis
functions ($n>m$) a WLS approximation is chosen as a
solution of equation~\eqref{eq:basic-system}, which becomes an overdetermined
problem. An example of this approach is DAM~\cite{Prax} that was originally
formulated to solve fluid flow in porous media. DAM uses six monomials for basis
and nine noded support domains to evaluate first and second derivatives of
physical fields required to solve the problem at hand, namely the Navier Stokes
equation. To determine the approximation coefficients $\b \alpha$, a norm
\begin{equation}
  R^2=\sum_j^n{w(\vp_j)(u(\vp_j)-\hat{u}(\vp_j))}^2 = (\b{B\alpha}-\b{u})^\T \b{W}^2 (\b{B\alpha} -\b{u}),
  \label{eq:residual-definition}
\end{equation}
is minimized, where $\b{W}$ is a diagonal matrix with elements $W_{jj}=\sqrt{w(\vp_j)}$ with
\begin{equation}
  w(\vec{p}) = \exp\left(-\left(\frac{\|\vec{p}_0-\vec{p}\|}{\sigma p_{\text{min}}}\right)^2\right),
  \label{eq:gaussian}
\end{equation}
where $\sigma$ stands for weight shape parameter, $\vec{p}_0$ for centre of
support domain  and $p_{\text{min}}$ for the distance to the nearest support
domain node.  There are different computational approaches to
minimizing~\eqref{eq:residual-definition}.  The most intuitive and also
computationally effective is to simply compute the gradient of $R^2$ with
respect to $\b{\alpha}$ and setting it to zero, resulting in a positive definite system
\begin{equation}
  \b{B}^\T \b W^2 \b B \b{\alpha} = \b{B}^\T \b W^2 \b{u}.
  \label{eq:wls-system}
\end{equation}
The problem of this approach is bad conditioning, as the condition
number of $\b{B}^\T \b W^2 \b B$ is the square of the condition number of $\b{WB}$,
unnecessarily increasing numerical instability.
A more stable and more expensive approach is QR decomposition. An even more
stable approach is SVD decomposition, which is of course even
more expensive.  Nevertheless, the solution of equation~\eqref{eq:wls-system}
can be written generally in matrix form as
\begin{equation}
  \b{\alpha} = (\b{W} \b{B})^+ \b{W} \b{u},
  \label{eq:alpha-def}
\end{equation}
where $\b A^+$ stands for a Moore–Penrose pseudo inverse of matrix $\b A$. By
explicitly inserting equation~\eqref{eq:alpha-def} for $\b{\alpha }$
into~\eqref{eq:u-definition}, the equation
\begin{equation}
  \label{eq:approx-system}
  \hat{u}(\vp) = \b{b}(\vp)^\T (\b{W} \b{B})^+\b{W}\b{u} = \b{\chi}(\vp)\b{u},
\end{equation}
is obtained, where $\b{\chi} = \b{b}(\vp)^\T (\b{W} \b{B})^+\b{W}$ is called a
shape function. Now, we can apply a partial differential operator $\L$ to the
trial function, and get
\begin{equation}
  (\L\hat{u})(\vp) = (\L\b{\chi})(\vp)\b{u}.
\end{equation}

In this paper we deal with a Cauchy-Navier equation and therefore following shape
functions are needed, expressed explicitly as
\begin{align}
  \label{eq:shapes-1}
  \b{\chi}^{\partial x}(\vp) &= \dpar{\b{b}}{x}(\vp)^\T (\b{W} \b{B})^+\b{W}, \\
  \b{\chi}^{\partial y}(\vp) &= \dpar{\b{b}}{y}(\vp)^\T (\b{W} \b{B})^+\b{W}, \\
  \b{\chi}^{\partial x \partial x}(\vp) &= \dpar{^2\b{b}}{x^2}(\vp)^\T (\b{W} \b{B})^+\b{W}, \\
  \b{\chi}^{\partial x \partial y}(\vp) &= \dpar{^2\b{b}}{x\partial y}(\vp)^\T(\b{W} \b{B})^+\b{W}, \\
  \b{\chi}^{\partial y \partial y}(\vp) &= \dpar{^2\b{b}}{^2y}(\vp)^\T(\b{W} \b{B})^+\b{W}.  \label{eq:shapes-5}
\end{align}
The shape functions depend only on the numerical setup, namely nodal
distribution, shape parameter, basis and support selection, and can as
such be precomputed for a specific computation.

\section{Governing problem}
\label{sec:governing-problem}
The goal in this paper is to numerically determine the stress and displacement
distributions in a solid body subjected to the applied external force. To obtain
a displacement vector field  $\vec{u}$ throughout the domain, a Cauchy-Navier
equation is solved, which can expressed concisely in vector form as
\begin{equation}
  \label{eq:navier}
  (\lambda + \mu) \nabla(\nabla \cdot  \vec{u}) + \mu \nabla^2 \vec{u} = 0,
\end{equation}
where $\mu$ and $\lambda$ stand for Lam\'e constants.
In two dimensions we express $\vec{u} = (u, v)$ and the equation
reads
\begin{align}
(\lambda+\mu) \frac{\partial}{\partial x}\left( \frac{\partial u}{\partial x} + \frac{\partial v}{\partial y} \right) + \mu \left( \frac{\partial^2 u}{\partial x^2} + \frac{\partial^2 u}{\partial y^2} \right) &= 0 \label{eq:navier-u}\\
(\lambda+\mu) \frac{\partial}{\partial y}\left( \frac{\partial u}{\partial x} + \frac{\partial v}{\partial y} \right) + \mu \left( \frac{\partial^2 v}{\partial x^2} + \frac{\partial^2 v}{\partial y^2} \right) &= 0 \label{eq:navier-v}
\end{align}
Two types of boundary conditions are commonly used when solving these types of problems, namely
essential boundary conditions and traction (also called natural) boundary
conditions. Essential boundary conditions specify displacements on some portion
of the boundary of the domain, i.e.~$\vec{u} = \vec{u}_0$, while traction
boundary conditions specify surface traction $\sigma\vec{n} = \vec{t}_0$, where
$\vec{n}$ is an outside unit normal to the boundary of the domain and
\begin{equation}
  \sigma = \begin{bmatrix} \sigma_{xx} & \sigma_{xy} \\
   \sigma_{xy} & \sigma_{yy} \end{bmatrix}
\end{equation}
is the stress tensor. In terms of
displacement vector $\vec{u}$ the traction boundary conditions read
\begin{align}
{t_0}_1 &= \mu n_2\dpar{u}{y} + \lambda n_1 \dpar{v}{y} + (2\mu+\lambda)n_1\dpar{u}{x} + \mu n_2 \dpar{v}{x} \label{eq:trac1} \\
{t_0}_2 &= \mu n_1\dpar{u}{y} + (2\mu+\lambda) n_2 \dpar{v}{y} + \lambda n_2\dpar{u}{x} + \mu n_1 \dpar{v}{x} \label{eq:trac2}
\end{align}
where ${t_0}_i$ and $n_i$
denote the Cartesian components of $\vec{t}_0$ and $\vec{n}$.

\section{Solution procedure}
\label{sec:solution-procedure}

\subsection{Discretization of the problem}
\label{sec:assembly}
The elliptic boundary value problem at hand is discretized into a linear system of
$2N$ algebraic equations by approximating the differential operations using MLSM,
as described in section~\ref{sec:mlsm}. A block system of linear equations for two
vectors $\b{u}$ and $\b{v}$ of unknowns representing values $u(\vec{p_i})$ and $v(\vec{p_i})$,
respectively, is constructed. This system is a discrete analogy of
PDE~(\ref{eq:navier}) and can symbolically be represented as
\begin{equation}
\label{eq:global-system}
\begin{bmatrix}
U1 & V1 \\ U2 & V2
\end{bmatrix}
\begin{bmatrix}
\b{u} \\ \b{v}
\end{bmatrix} =
\begin{bmatrix}
\b{b_1} \\ \b{b_2}
\end{bmatrix},
\end{equation}
where $\b{u}$ and $\b{v}$ stand for unknown displacements,
$\b{b_1}$ and $\b{b_2}$ for values of
boundary conditions and blocks $U1, V1, U2, V2$ contain precomputed shape
functions~(\ref{eq:shapes-1}--\ref{eq:shapes-5}).
With $\mathcal{N}(i)$ standing for a list of indices of the chosen $n$
neighbours of a point $\vec{p}_i$, as introduced in the beginning of
section~\ref{sec:mlsm}, we can, for all indices $i$ of internal nodes, express
\begin{eqnarray}
\left.
\begin{aligned}
U1_{i, \mathcal{N}(i)_j} &= \left[ (\lambda + 2\mu)\b{\chi}^{\partial x\partial x}(\vec{p}_i) + \mu\b{\chi}^{\partial y\partial y}(\vec{p}_i) \right]_j \label{eq:u1} \\
V1_{i, \mathcal{N}(i)_j} &= \left[ (\lambda + \mu)\b{\chi}^{\partial x\partial y}(\vec{p}_i) \right]_j  \\
{\b{b}_1}_{i} &= 0
\end{aligned}
\right\},
\end{eqnarray}
\begin{eqnarray}
\left.
\begin{aligned}
U2_{i, \mathcal{N}(i)_j} &= \left[ (\lambda + \mu)\b{\chi}^{\partial x\partial y}(\vec{p}_i) \right]_j \label{eq:u2}  \\
V2_{i, \mathcal{N}(i)_j} &= \left[ \mu\b{\chi}^{\partial x\partial x}(\vec{p}_i) + (\lambda + 2\mu)\b{\chi}^{\partial y\partial y}(\vec{p}_i) \right]_j \\
{\b{b}_2}_{i} &= 0 \\
\end{aligned}
\right\},
\end{eqnarray}
for each $j = 1, \ldots, n$. Note that equation~\eqref{eq:u1}  represents direct discrete analogue of~\eqref{eq:navier-u} and, likewise, \eqref{eq:u2} of~\eqref{eq:navier-v}.

Similarly, for all indices $i$ of boundary nodes with traction boundary conditions we express
\begin{eqnarray}
\left.
\begin{aligned}
  U1_{i, \mathcal{N}(i)_j} &= \left[ \mu n_2 \b{\chi}^{\partial y}(\vec{p}_i) + (2\mu + \lambda)n_1 \b{\chi}^{\partial x}(\vec{p}_i) \right]_j \label{eq:t1}  \\
V1_{i, \mathcal{N}(i)_j} &= \left[ \lambda n_1\b{\chi}^{\partial y}(\vec{p}_i) +  \mu n_2 \b{\chi}^{\partial x}(\vec{p}_i) \right]_j  \\
{\b{b}_1}_{i} &= {t_0}(\vec{p}_i)_1
\end{aligned}
\right\},
\end{eqnarray}
\begin{eqnarray}
\left.
\begin{aligned}
  U2_{i, \mathcal{N}(i)_j} &= \left[ \mu n_1 \b{\chi}^{\partial y}(\vec{p}_i) + \lambda n_2 \b{\chi}^{\partial x}(\vec{p}_i) \right]_j \label{eq:t2}  \\
V2_{i, \mathcal{N}(i)_j} &= \left[ \mu n_1\b{\chi}^{\partial x}(\vec{p}_i) + (2\mu + \lambda)n_2 \b{\chi}^{\partial y}(\vec{p}_i) \right]_j  \\
{\b{b}_2}_{i} &= {t_0}(\vec{p}_i)_2
\end{aligned}
\right\},
\end{eqnarray}
for each $j = 1, \ldots, n$, where $n_i$ are the Cartesian components of
the outside unit normal to the boundary in node $\vec{p}_i$.
Again, equation~\eqref{eq:t1} is a direct analogue of~\eqref{eq:trac1}
and~\eqref{eq:t2} of~\eqref{eq:trac2}.
And finally, for indices $i$ of nodes with essential boundary condition, we express
\begin{eqnarray}
\begin{aligned}
U1_{i, i} &=1  \\
{\b{b}_1}_{i} &= {u_0}(\vec{p}_i)_1
\end{aligned}
& \quad \text{ and } \quad  &
\begin{aligned}
U2_{i, i} &=1  \\
{\b{b}_2}_{i} &= {u_0}(\vec{p}_i)_2
\end{aligned}
\quad .
\end{eqnarray}

\begin{figure}[h]
  \centering
  \includegraphics[width=0.5\textwidth]{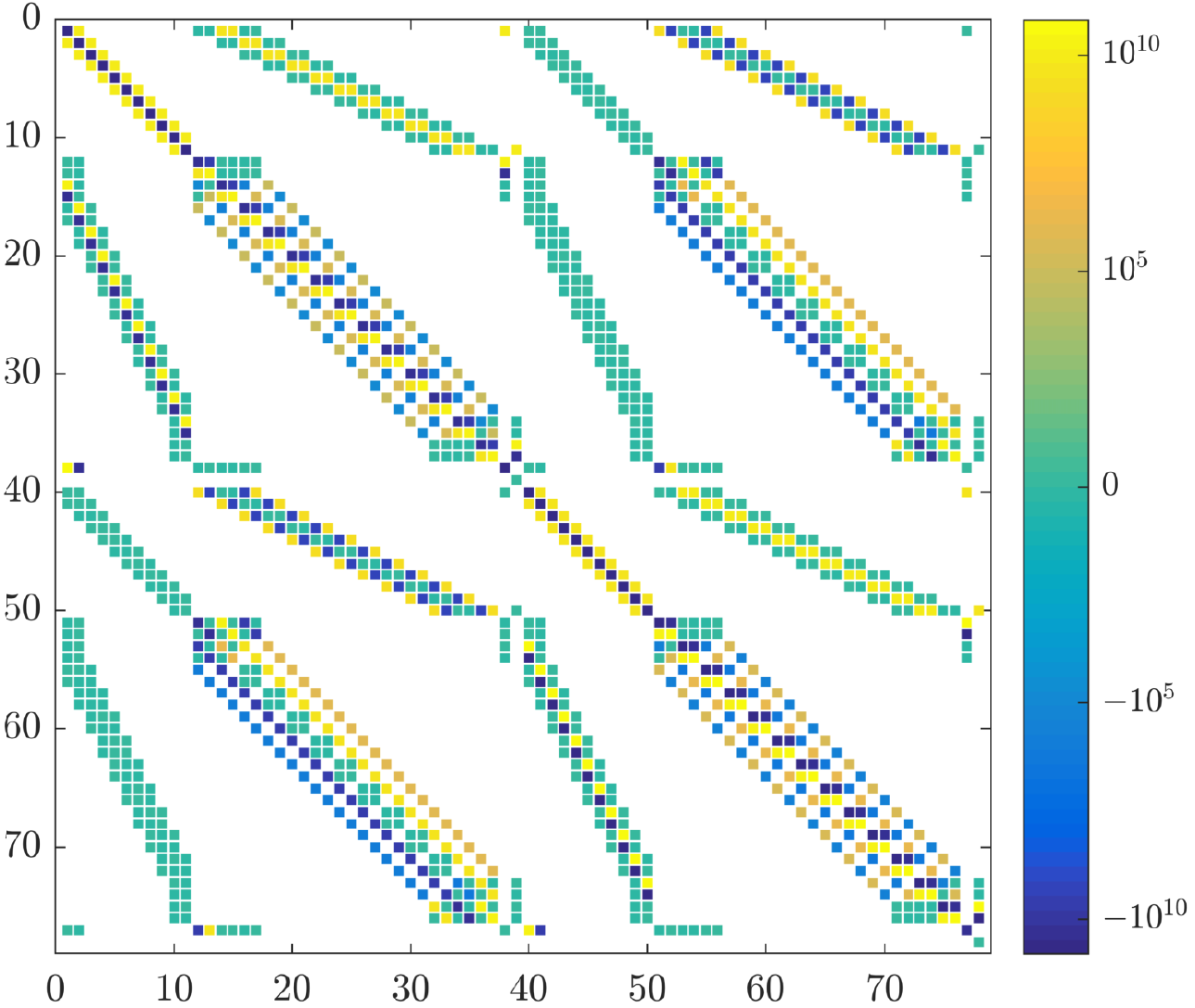}
  \caption{Matrix of the final system of equations in cantilever beam case with
    $N = 39$ and \unit[22]{\%} nonzero elements.}
  \label{fig:matrix}
\end{figure}

System~\eqref{eq:global-system} is sparse with nonzero ratio of
less then $2n/N$. An example of the matrix of this system for the
cantilever beam problem described in section~\ref{sec:cantilever-beam}
is shown in Figure~\ref{fig:matrix}, where the block structure
and different patterns for boundary and internal nodes are clearly visible.

\subsection{Positioning of nodes in a complex domain}
\label{sec:positioning}
Meshless methods are advertised as the methods that do not require any
topological relations among nodes. That implies that even randomly distributed
nodes could be used~\cite{random_mesh}.  However, it is well-known that with
regularly distributed nodes one achieves much better results in terms of
accuracy and stability~\cite{amani}. This has also been recently reported for
MLSM in a solution of a Navier-Stokes problem~\cite{Kouszeg1}. The reason behind
the sensitivity regarding the distribution of nodes lies in the generation of shape
functions. To construct a stable method well balanced support domains are
needed, i.e.\ the nodes in support domain need to be distributed evenly enough
\cite{Kouszeg1}. This condition is obviously fulfilled in regular nodal
distributions, but when working with more interesting geometries, the
positioning of nodes requires additional treatment. In literature one can find
several algorithms for distributing the nodes within the domain of different
shapes~\cite{onate, bubble}. In this paper we will use an extremely simple
algorithm, introduced in~\cite{Kouszeg1} to minimize the variations in distances
between nodes in the support domain. The basic idea is to ``relax'' the nodes
based on a potential between them. Since a Gaussian function is a suitable
potential and already used as weight in the shape functions, the nodes are
translated simply as
\begin{equation}
  \delta  \vec{p}(\vec{p}) = -\sigma_k \sum_{i=1}^{N_S} \nabla w(\vec{p}-\vec{p_i}),
\end{equation}
where $\delta \vec{p}$, $\vec{p}_i$, $\sigma_k$ and $N_S$ stand for the
translation step of the node, position of $i$-th support node, relaxation
parameter and number of support nodes, respectively (Figure
\ref{fig:h-relax-alg}). After offsets in all nodes are computed, the nodes are
repositioned as
\begin{equation}
  \vec{p} \gets \vec{p} + \delta \vec{p}(\vec{p}).
\end{equation}
Presented iterative process procedure begins by positioning the boundary nodes,
which are considered as the definition of the domain and are kept fixed throughout
the process.

\subsection{h-refinement}
\label{sec:h-ref}
Besides flexibility regarding the shape of the domain, nodal refinement is often
mandatory to achieve desired accuracy in cases with pronounced differences in
stress within the domain. A typical example of such situation is a contact
problem~\cite{pereira2016convergence}. To mitigate the error in areas with high
stresses the h-refinement scheme is used., which has already
been introduced into different meshless solutions~\cite{libre, Bourantas}. In
context of local RBF approximation the h-refinement has been used in the solution of
the Burger's equation~\cite{Kouszeg_cmc}, where a quad-tree based algorithm has
been used to add and remove child nodes symmetrically around the parent node in
transient solution of Burgers' equation. However, the algorithm presented in ~\cite{Kouszeg_cmc} supported only regular nodal distribution. In this paper we generalize it also to irregular nodal distribution.

\begin{figure}[h!]
  \centering
  \begin{subfigure}[t]{0.45\textwidth}
    \includegraphics[width=1\textwidth]{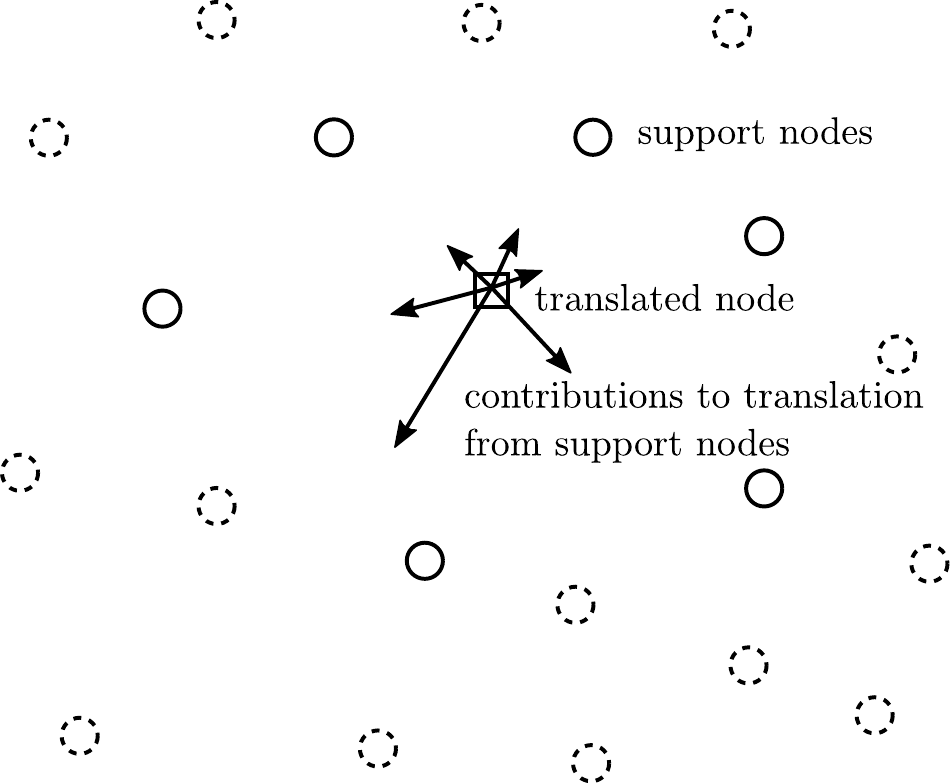}
    \caption{Scheme of the relax algorithm.}
    \label{fig:h-relax-alg}
  \end{subfigure}
  \begin{subfigure}[t]{0.45\textwidth}
    \includegraphics[width=1\textwidth]{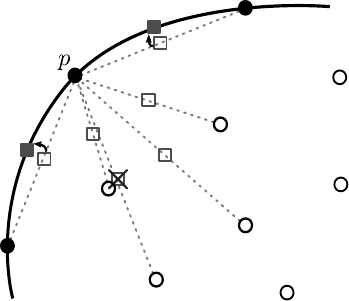}
    \caption{Scheme of the h-refinement algorithm.}
    \label{fig:h-refine-alg}
  \end{subfigure}\hspace{5mm}
  \caption{Schemes of algorithms used to improve the quality of the discretization.}
\end{figure}

In each node to be refined, new nodes are added on the half
distances between the node itself and its support nodes
\begin{equation}
 \vec{p}_{j}^{\;new} = \frac{\vec{p} + \vec{p}_j}{2},
\end{equation}
 where index $j$ indicates $j$-th support node. When
adding new nodes, checks are performed if the newly added node is too close to any
of the existing nodes; in that case the node is not added. Moreover, if the refined
node and support node are both boundary nodes, newly added node is positioned on
the boundary (Figure~\ref{fig:h-refine-alg}). This procedure can be repeated several
times if an even more refined domain is desired. These subsequent refinements will
be called levels of refinement and will be denoted as level $i$ for the
refinement that resulted from $i$ applications of the described algorithm.

The described algorithm follows the concept
of meshless methods and as such does not require any special topological
relations between nodes to refine a certain part of the computation domain. It
is also flexible regarding the dimensionality of the domain, i.e.\ there is no
difference in implementation of 2D or 3D variant of the algorithm.

\begin{figure}[h!]
  \centering
  \includegraphics[width=0.7\textwidth]{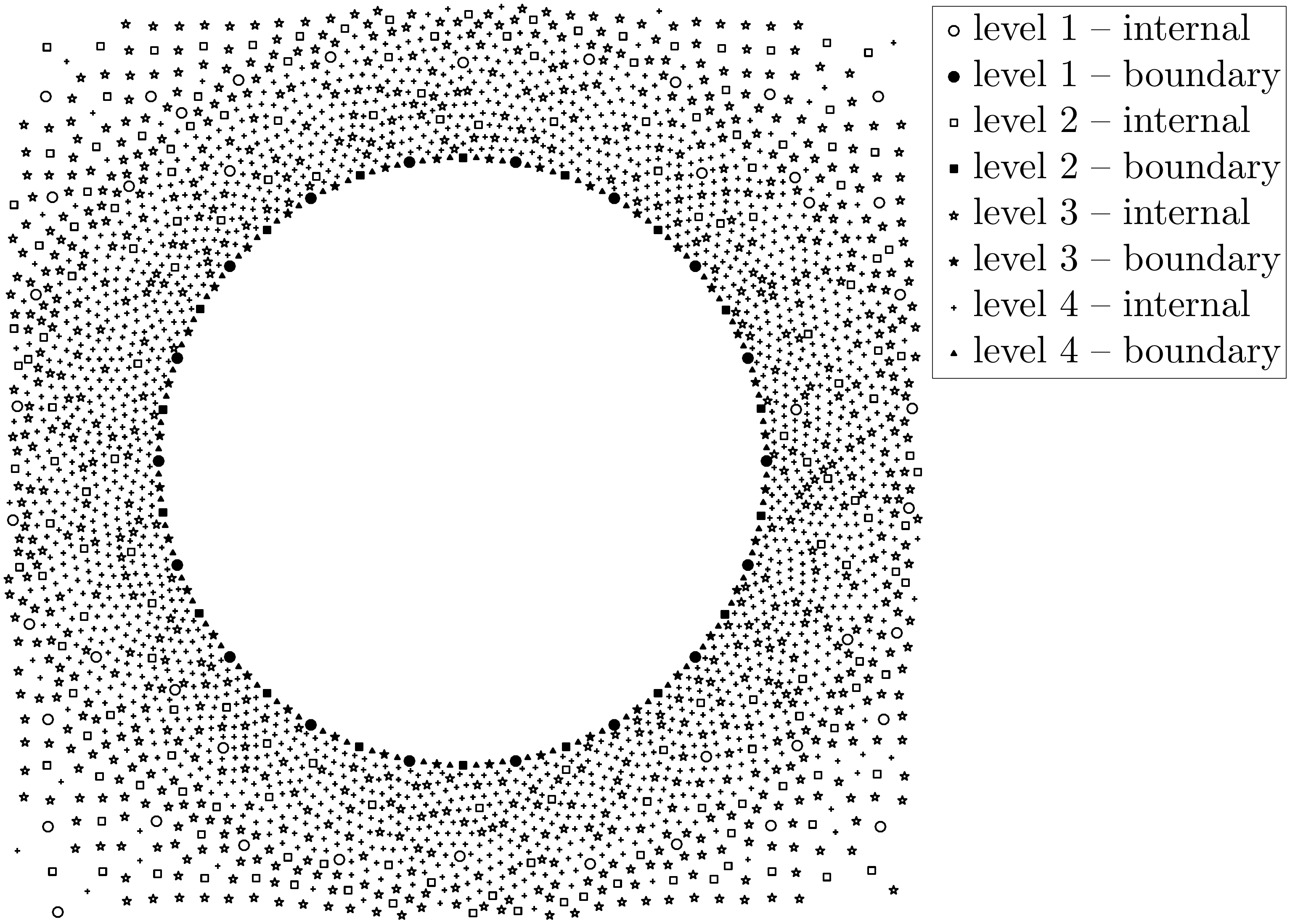}
  \caption{Four levels of the refinement algorithm applied around a hole in a domain after relaxation.}
  \label{fig:refine-hole}
\end{figure}

An example on a non-trivial refinement is demonstrated in Fig.~\ref{fig:refine-hole},
where a domain with a hole is considered. The vicinity of the hole is four times
refined and then, to mitigate possible irregularities during refinement, relaxed.

\subsection{Asymptotic complexity of MLSM}
\label{sec:complexity}
The asymptotic complexity analysis begins with an assumption that evaluations of
basis functions, weights, linear operators and boundary conditions take $O(1)$
time. For simple domain discretization, such as the uniform grid in a rectangle
or random positioning, $O(N)$ time is required, where $N$ stands for number of
computational nodes. To find the neighbours of each
point, a tree based data structure such as kd-tree~\cite{trobec_kosec}, taking
$O(N \log N)$ time to construct and $O(n \log N)$ time to query $n$ closest
nodes, is used. The relaxation of nodal positions (see
section~\ref{sec:positioning}) with $I$ iterations costs additional $O(I n N
\log^2 N)$ time. Re-finding the support nodes by rebuilding
the tree and querying for support nodes once again, requires another $O((N+n)\log N)$
time. Calculation of the shape functions requires $N$ SVD decompositions, each
taking $O(nm^2)$ time, as well as some matrix and vector multiplication of lower
complexity. Assembling the matrix takes $O(nN)$ time and assembling the right
hand side takes $O(N)$ of time. Then, the system is solved using BiCGSTAB iterative
algorithm.  The final time complexity is thus $O(InN \log^2 N + (N+n)\log N + m^2n N)
+ T$, where $T$ stands for the time spent by BiCGSTAB.


For comparison, the complexity of a well-known weak form Element Free Galerking
method (EFG)~\cite{EFG} differs from MLSM in construction of the shape functions,
whose computation requires $O(N n_q m^2 n)$ time using EFG method,
with $n_q$ standing for the number of Gauss integration points per node.
Additionally, the number of nonzero elements in the final system of EFG is
of order $n_q$ times higher than that of MLSM, again increasing the complexity of EFG.

\section{Numerical examples}
\label{sec:results}
 \subsection{Cantilever beam}
\label{sec:cantilever-beam}
First, the standard cantilever beam test is solved to assess accuracy and
stability of the method.  Consider an ideal thin cantilever beam of length $L$
and height $D$ covering the area $[0, L] \times [-D/2, D/2]$.  Timoshenko beam
theory offers a closed form solution for displacements and stresses in such a
beam under plane stress conditions and a parabolic load on the left side. The
solution is widely known and derived in e.g.~\cite{slaughter2012linearized},
giving stresses in the beam as
\begin{align}
\sigma_{xx} &= \frac{Pxy}{I}, \quad \sigma_{yy} = 0, \quad \sigma_{xy} = \frac{P}{2I} \left( \frac{D^2}{4}
- y^2 \right),
\label{eq:cantilever-beam-analytical-stress}
\end{align}
and displacements as
\begin{align}
u &= \frac{P y \left(3 D^2 (\nu +1)-4 \left(3 L^2+(\nu +2) y^2-3 x^2\right)\right)}{24 E I},
\label{eq:cantilever-beam-analytical} \\
v &= -\frac{P \left(3 D^2 (\nu +1) (L-x)+4 (L-x)^2 (2 L+x)+12 \nu  x y^2\right)}{24 E I},
\nonumber
\end{align}
where $I = \frac{1}{12} D^3$ is the moment of inertia around the horizontal
axis, $E$ is Young's modulus, $\nu$ is the Poisson's ratio and $P$ is the
total load force.

In the numerical solution, traction free boundary conditions are used on
the top and bottom of the domain, essential boundary conditions given
by~\eqref{eq:cantilever-beam-analytical} are used on the right and traction
boundary conditions given by~\eqref{eq:cantilever-beam-analytical-stress}
on the left
\begin{align}
  u(L, y) &= \frac{Py (2D^2(1+\nu) - 4(2+\nu) y^2)}{24 E I} \\
  v(L, y) &= -\frac{L \nu P y^2}{2 E I}
\end{align}
\begin{align}
  \mu\dpar{u}{y}(x, D/2) + \mu\dpar{v}{x}(x, D/2)) &= 0 \\
  \lambda \dpar{u}{x}(x, D/2) + (\lambda + 2\mu) \dpar{v}{y}(x, D/2) &= 0 \\
  -\mu\dpar{u}{y}(x, -D/2) - \mu\dpar{v}{x}(x, -D/2)) &= 0 \\
   -\lambda \dpar{u}{x}(x, -D/2) - (\lambda + 2\mu) \dpar{v}{y}(x, -D/2) &= 0 \\
  -\lambda \dpar{v}{y}(0, y) - (\lambda + 2\mu) \dpar{u}{x}(0, y)&= 0 \\
  -\mu\dpar{u}{y}(0, y) - \mu\dpar{v}{x}(0, y) &= \frac{P}{2 I} ((D/2)^2-y^2).
\end{align}
The problem is solved using MLSM method with $n = 9$ or $n = 13$ support nodes
and Gaussian weight with $\sigma = 1$ (see \eqref{eq:gaussian}).
Two sets of basis functions are considered, 9 monomials
\begin{equation}
  \b{b} =\{1, x, y,  x^2, y^2, xy, x^2y, xy^2, x^2y^2\}
\end{equation}
and 9 Gaussian basis functions (see~\eqref{eq:gaussian} for definition)
centred in support nodes.
In the following discussions these two choices of basis functions
will be referred to as M9 and G9, respectively.

System~\eqref{eq:global-system} is solved with BiCGSTAB iterative
algorithm~\cite{van1992bi} with ILUT preconditioning~\cite{saad1994ilut}.
Values of $L = \unit[30]{m}$, $D = \unit[5]{m}$, $E = \unit[72.1]{GPa}$,
$\nu = 0.33$ and $P = \unit[1000]{N/m}$ were chosen as physical parameters
of the problem.

\begin{figure}[!ht]
  \centering
  \includegraphics[width=0.7\textwidth]{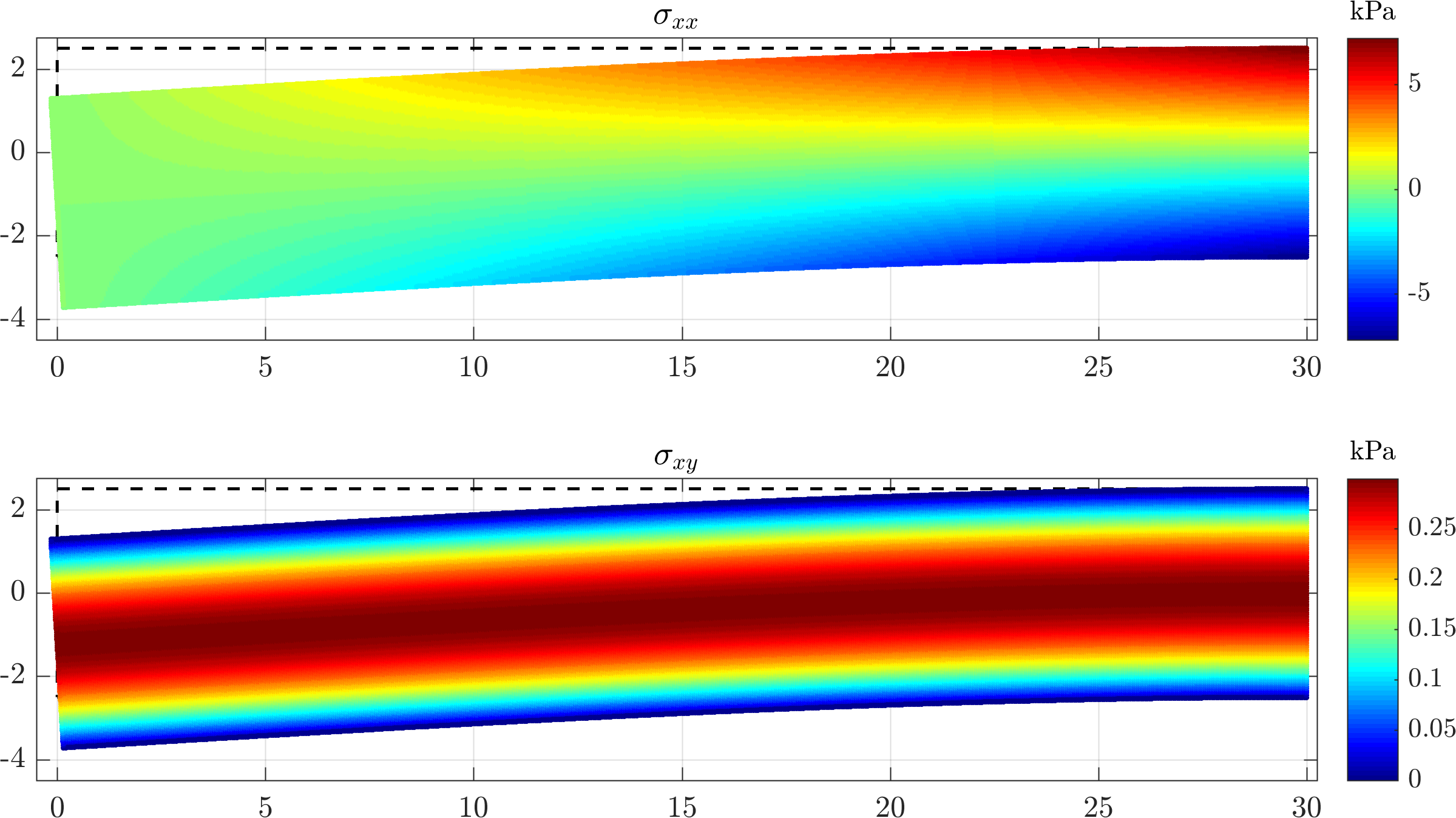}
  \caption{Numerical solution of cantilever beam case. Note that for the sake of
  visibility the displacements are multiplied by factor $10^5$.}
  \label{fig:cantilever-beam-solution}
\end{figure}

The acquired numerical solution of the cantilever beam problem is shown in
Figure~\ref{fig:cantilever-beam-solution}.

\begin{figure}[!ht]
  \centering
  \includegraphics[width=0.49\textwidth]{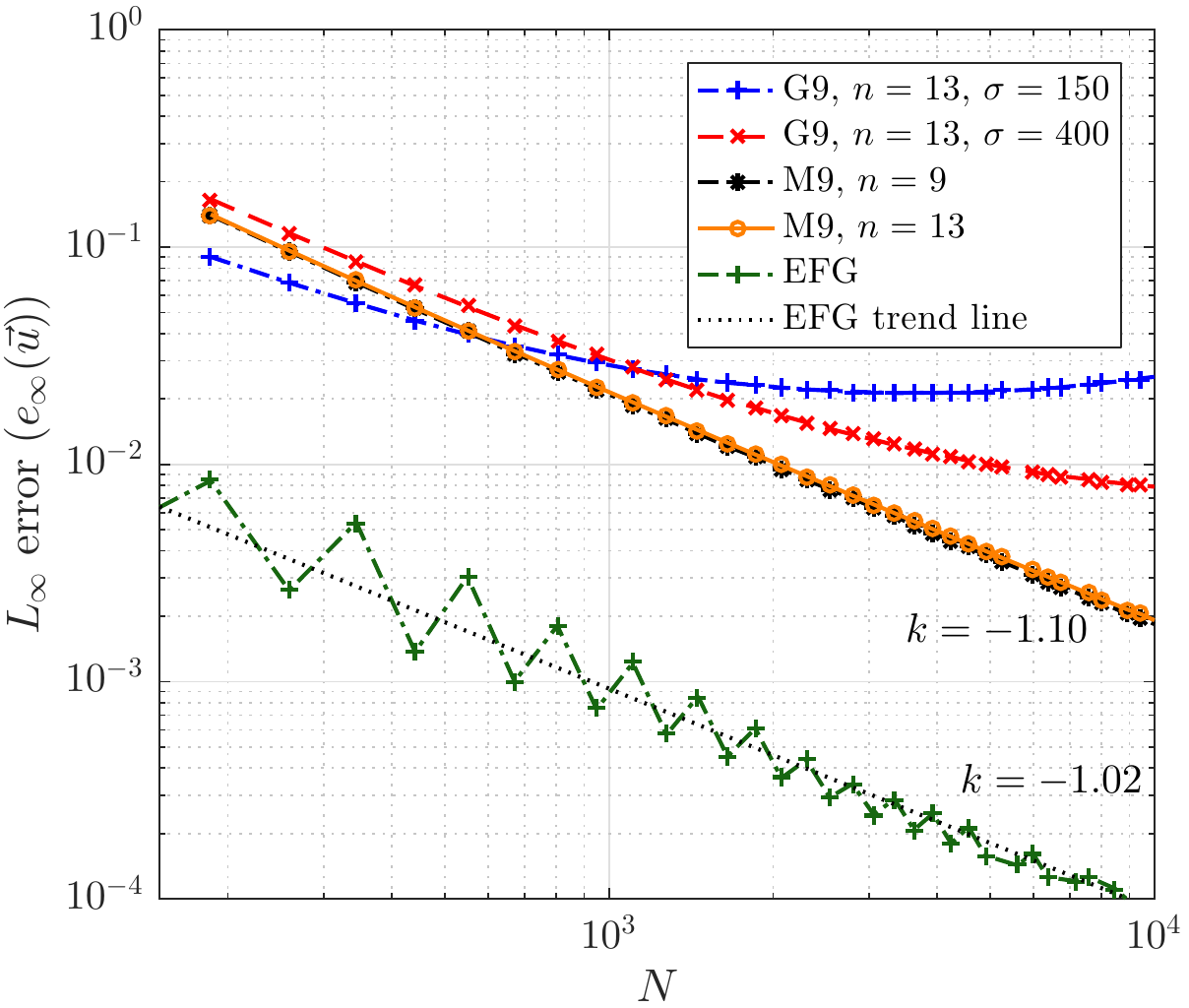}
  \includegraphics[width=0.49\textwidth]{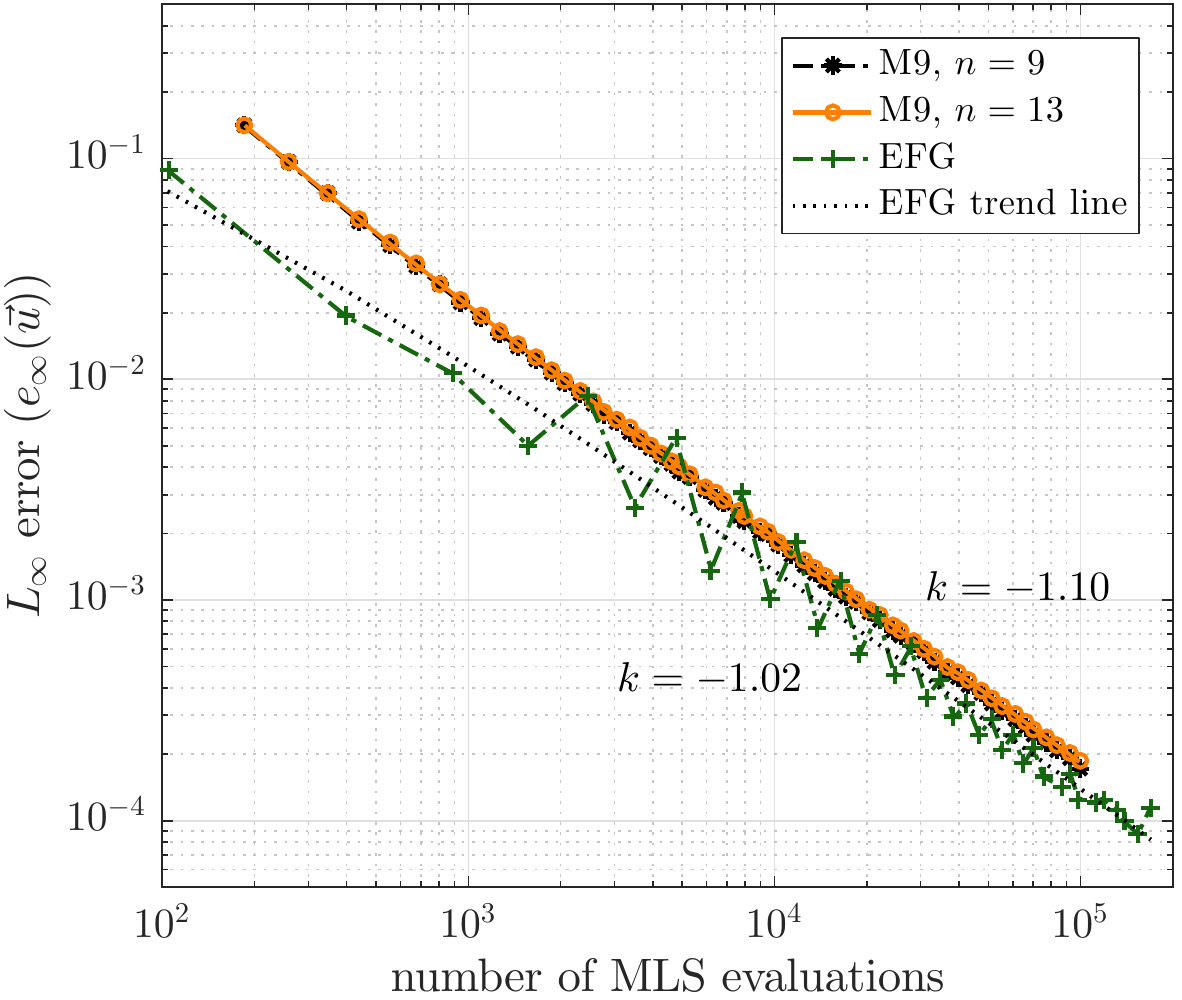}
  \caption{Accuracy of different MLSM setups compared to EFG per number of
    computational nodes (left) and number of MLS evaluations (right).}
  \label{fig:cantilever-beam-conv}
\end{figure}

The error of the numerical approximation of stresses and
displacements is measured in relative discrete $L_\infty$ norm, using
\begin{align}
  e_\infty(\vec{u}) &=
   \frac{\max_{x\in X} \{ \max\{|u(x)-\hat u(x)|,|v(x)-\hat v(x)|\}\}}
   {\max_{x\in X} \{\max\{|u(x)|,|v(x)|\}\}} \qquad \text{and} \\
  e_\infty(\sigma) &= \frac{\max_{x\in X} \{\max\{|\sigma_{xx}(x)-\hat{\sigma}_{xx}(x)|,
  |\sigma_{yy}(x)-\hat{\sigma}_{yy}(x)|,
  |\sigma_{xy}(x)-\hat{\sigma}_{xy}(x)| \}\}}{
  \max_{x\in X} \{\max\{|\sigma_{xx}(x)|, |\sigma_{yy}(x)|, |\sigma_{xy}(x)| \}\}},
\end{align}
as error indicators, with  $X$ representing the set of all nodes. Convergence
with respect to the number of computational nodes is shown in
Figure~\ref{fig:cantilever-beam-conv}. The numerical approximations converge
towards the correct solution in stress $(e_\infty(\sigma))$ norm as well,
with approximately the same convergence rate.

It can be seen that monomials converge very regularly with order 1 as expected,
while Gaussian functions exhibit slightly worse convergence. Such behaviour has
already been reported in solution of diffusion equation, where MLSM with Gaussian
basis failed to obtain accurate solution with a high number of computational nodes.
More details about the phenomenon and further reading can be found
in~\cite{KosecZinterhof}.

The method was compared to the standard Element Free Galerkin (EFG)
method~\cite{MMrev}. The EFG method used circular domains of influence with
radius $d_I$ equal to 3.5 times internodal distance, a cubic spline
\begin{equation}
  w(\vec{p}) = \tilde{w}\left(\frac{\|\vec{p} - \vec{p}_i\|}{d_I}\right), \quad
  \tilde{w}(r) =
  \begin{cases}
    \frac23 - 4r^2 + 4r^3 & 0 \leq r < \frac12 \\
    \frac43 - 4r + 4r^2 - \frac43 r^3 & \frac12 \leq r < 1 \\
    0 & 1 \leq r \\
  \end{cases}
\end{equation}
for a weight function, $n_q = 4$ Gaussian points for approximation of line integrals
and $n_q = 16$ points for approximating area integrals.
Lagrange multipliers were used to impose essential boundary conditions.

The performance of EFG with respect to the number of nodes is much better than
MLSM. However, a more fair comparison would take into account also a higher
complexity of the EFG. This can be achieved by comparing error with respect to
the number of MLS evaluations, which is the most time consuming part of the
solution procedure. In Figure~\ref{fig:cantilever-beam-conv} it is
demonstrated that although EFG provides much better results in comparison
to MLSM at a given number of nodes, its accuracy becomes comparable to MLSM,
when compared per number of MLS evaluations.

\begin{figure}[!ht]
  \centering
  \includegraphics[width=\textwidth]{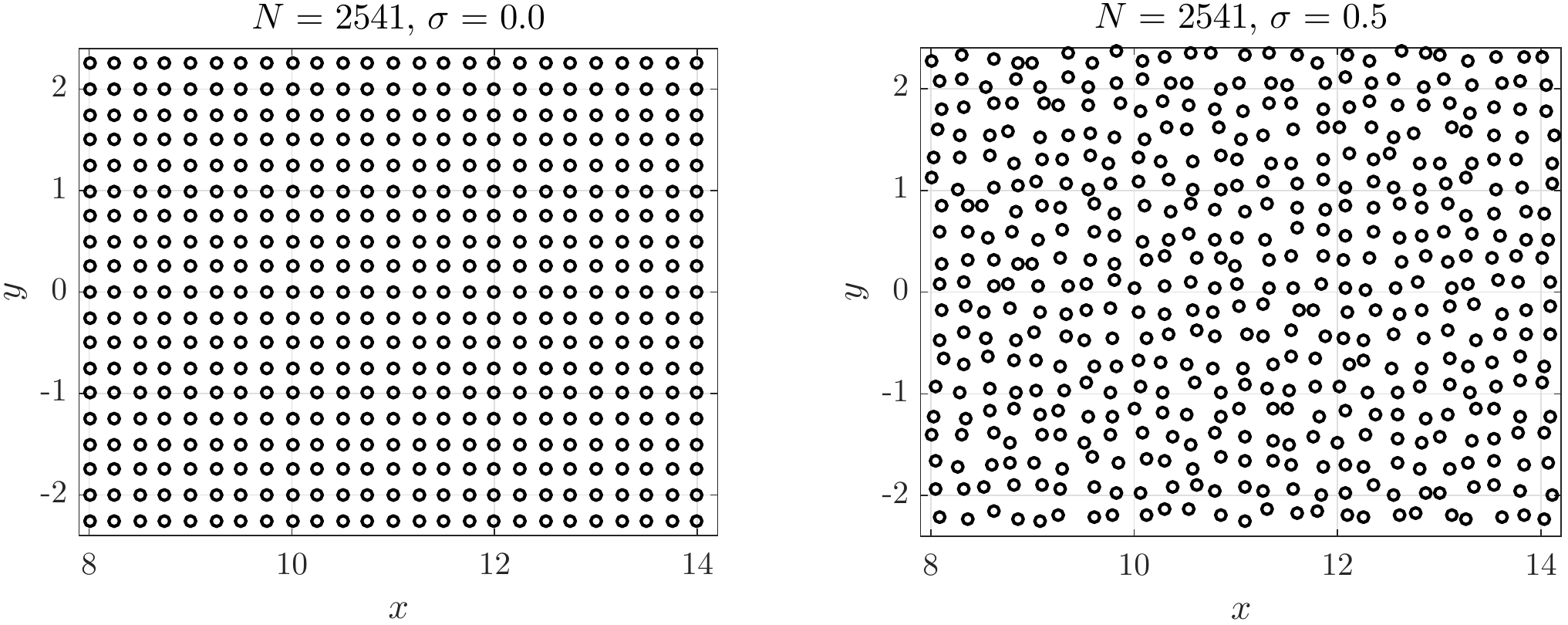}
  \caption{Regular and perturbed node positions, as used in stability analysis.}
  \label{fig:randomized-nodes}
\end{figure}

To asses the stability of the method regarding the nodal distribution,
the following analysis was performed. A regular distribution of points
as used in the solution in Figure~\ref{fig:cantilever-beam-solution} was distorted
by adding a random perturbation to each internal node. Its position is altered by
\begin{equation}
  \hat{\vec{p}} \gets \vec{p} + \sigma \vec{U}, \qquad \vec{U} \sim \text{Uniform}([0, \delta]^2),
\end{equation}
where $\delta$ is the distance to the closest node, and measuring the accuracy of
the solution with respect to $\sigma$, representing magnitude of the perturbation.
An example of original and perturbed node distributions are shown in
Figure~\ref{fig:randomized-nodes}.

\begin{figure}[!ht]
  \centering
  \includegraphics[width=0.7\textwidth]{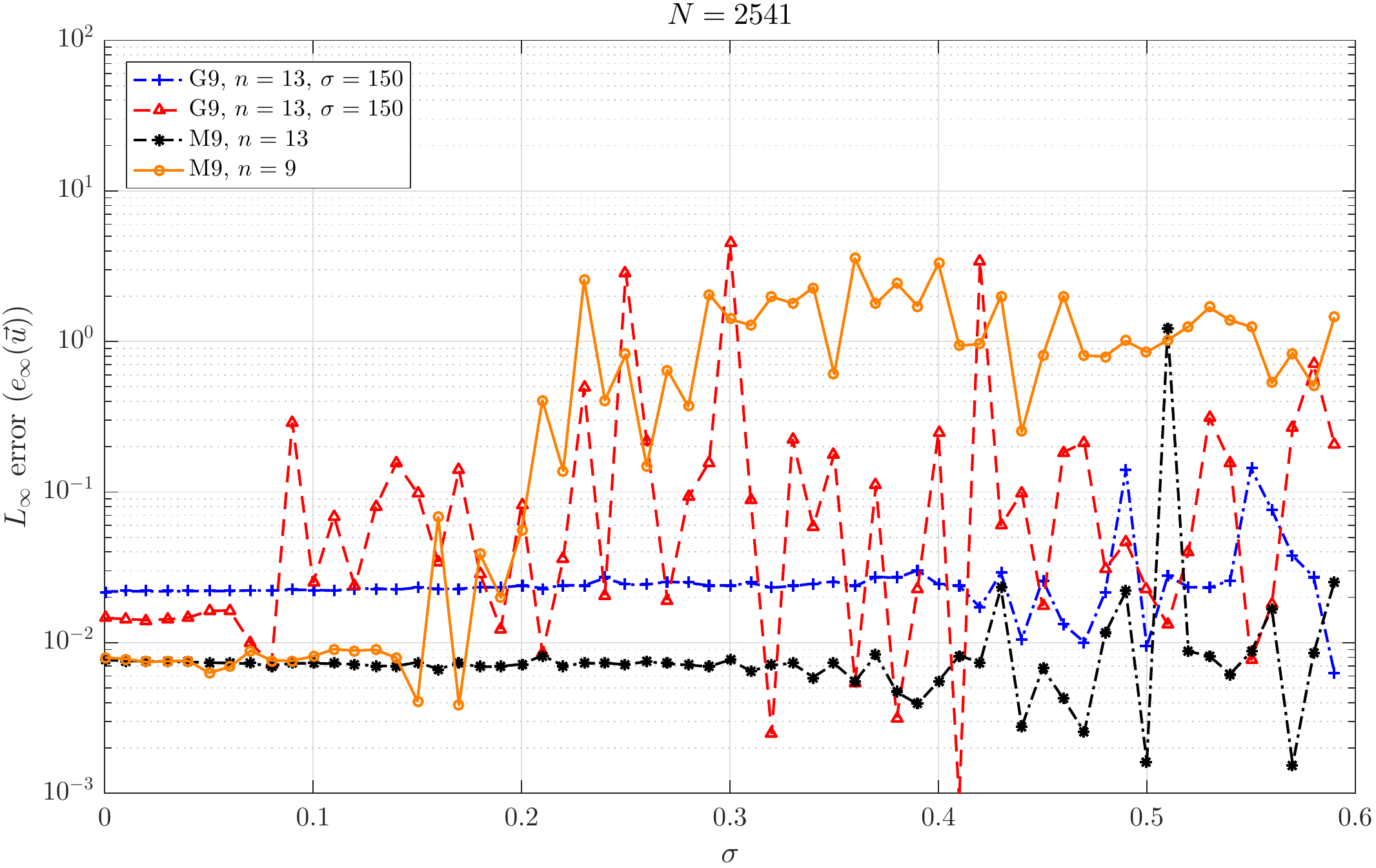}
  \caption{Stability of MLSM with respect to nodal perturbations.}
  \label{fig:cantilever-beam-stab}
\end{figure}

Accuracy of the solution with respect to the perturbation magnitude is presented in Figure~\ref{fig:cantilever-beam-stab}.
It is demonstrated that using monomials as a basis with 9 support nodes
results in an unstable setup. On the other hand monomials with 13 support nodes
are much more stable and equally accurate, while using Gaussian basis with high
shape parameter is the most unstable setup.
To mitigate the stability issue, a lower shape
parameter can be chosen, however, at the cost of accuracy. Regardless of the
setup one can expect the solution to be stable at least up to $\sigma \approx 0.1$.
Note that using more nodes in support domain can also increase stability. Refer to
\cite{Kouszeg1} for more details.

\begin{figure}[!ht]
  \centering
  \includegraphics[width=0.7\textwidth]{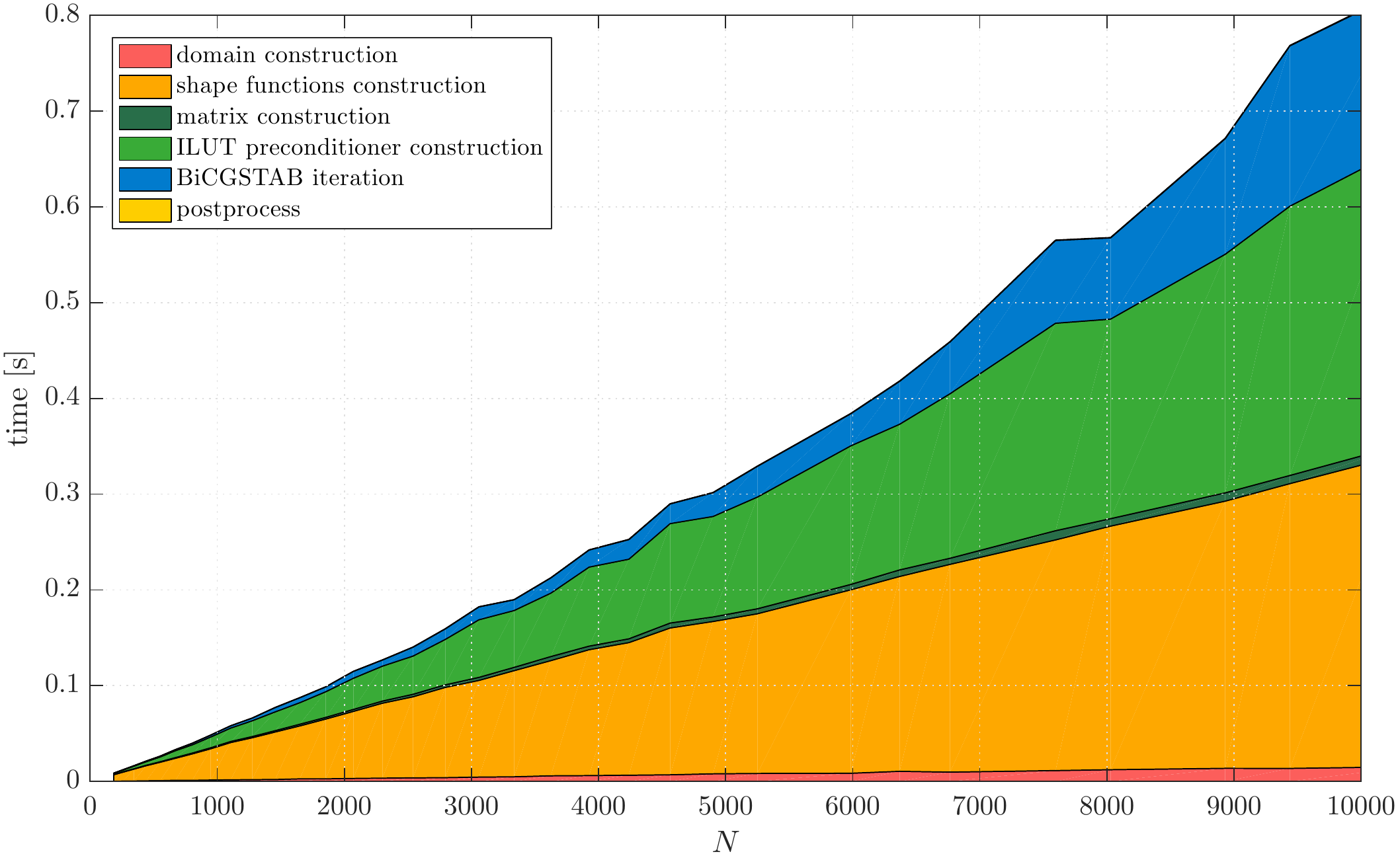}
  \caption{Execution time for different parts of the solution procedure with
    respect to the number of computational nodes.}
  \label{fig:cantilever-time}
\end{figure}

Time spent on each part of the solution procedure is shown in
Figure~\ref{fig:cantilever-time}. All measurements were performed on a laptop
computer with an \texttt{Intel(R) Core(TM) i7-4700MQ @2.40GHz} CPU and with
\unit[16]{GiB} of DDR3 RAM. MLSM is implemented in C++~\cite{utils_web}
and compiled using \verb|g++ 7.1.7| for Linux with
\texttt{-std=c++14 -O3 -DNDEBUG} flags.
It can be seen that solving the system~\eqref{eq:global-system} makes
up for more than \unit[50]{\%} of total time spent. Around \unit[70]{\%}
of that time is spent on computing the preconditioner. The only other
significant factor is computing the shape functions taking approximately
\unit[40]{\%} of total time. Domain construction and matrix assembly take
negligible amounts of time, matching the predictions made by
complexity analysis in section~\ref{sec:complexity}.

\begin{figure}[!ht]
  \centering
  \includegraphics[width=0.7\textwidth]{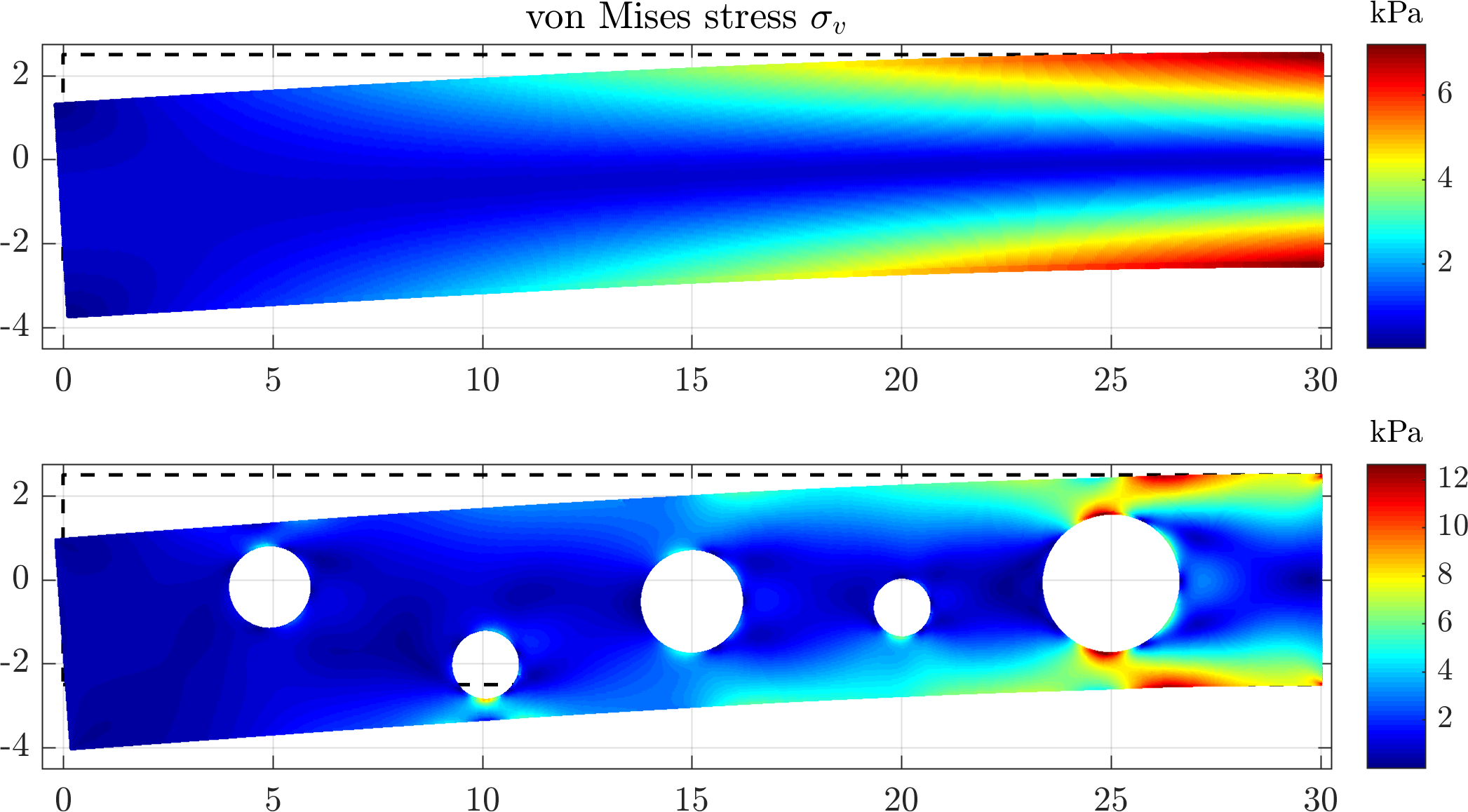}
  \caption{Numerical solution of a drilled cantilever beam case using
  $N = 177618$ nodes. Note that for the sake of visibility the displacements
  are multiplied by factor $10^5$.}
  \label{fig:cantilever-beam-holes}
\end{figure}

To emphasize the generality of MLSM method, a ``drilled'' domain is considered
in the next step. Arbitrarily positioned holes are added to the rectangular
domain. The positioning algorithm described in section~\ref{sec:positioning} and
h-refinement algorithm described in section~\ref{sec:h-ref} are used to
distribute the nodes inside the domain and refine the areas around the holes.
The boundary conditions in this example are  $\vec{u} = 0$ on the right,
traction free on the inside of the holes and on top and bottom and uniform load
of $P/D$ on the left. The computed solution is shown in
Figure~\ref{fig:cantilever-beam-holes} along with the ordinary cantilever beam
example.
Both solutions are coloured using von~Mises stress $\sigma_v$,  computed for the plane stress case as
\begin{equation}
  \sigma_v = \sqrt{\sigma_{xx}^2-\sigma_{xx}\sigma_{yy}+\sigma_{yy}^2+3\sigma_{xy}^2}.
\end{equation}

\begin{figure}[!ht]
  \centering
  \includegraphics[width=0.8\textwidth]{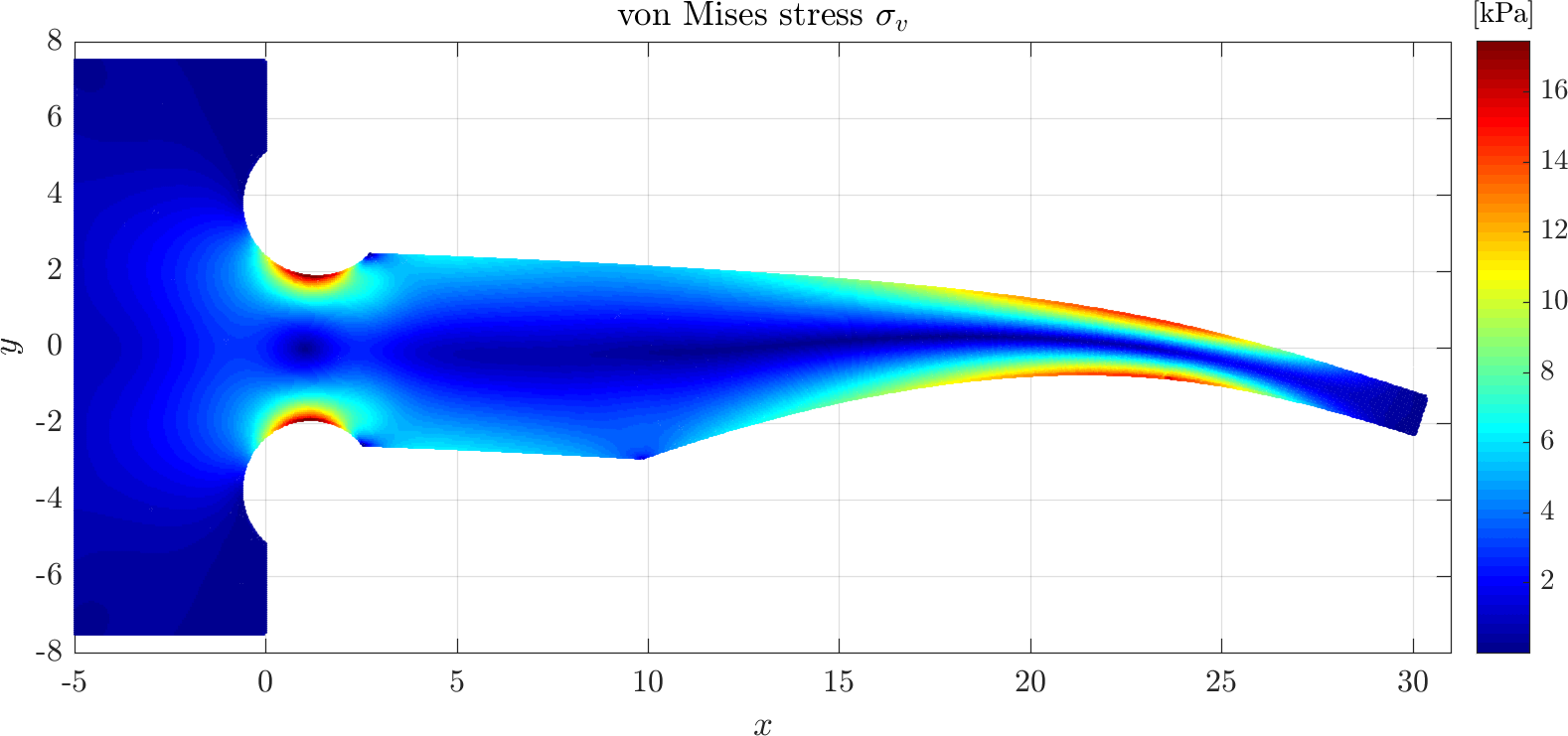}
  \caption{Numerical solution of a irregular cantilever beam using
    $N = 67887$ nodes. Note that for the sake of visibility the displacements
    are multiplied by factor $10^5$.}
  \label{fig:cantilever-beam-weird}
\end{figure}

To further illustrate the generality of the method, an even more deformed domain is
considered (Figure \ref{fig:cantilever-beam-weird}).


\subsection{Hertzian contact}
\label{sec:hertzian}
Another more interesting case arises from basic theory of contact mechanics, called
Hertzian contact theory~\cite{williams2001contact}.  Consider
two cylinders with radii $R_1$ and $R_2$ and parallel axes pressed together by a
force per unit length of magnitude $P$. The theory predicts they form a small
contact surface of width $2b$, where
\begin{equation}
  b = 2\sqrt{\frac{PR}{\pi E^*}}
\end{equation}
and
\begin{align}
  \frac{1}{R} &= \frac{1}{R_1} + \frac{1}{R_2}, \\
  \frac{1}{E^*} &= \frac{1-{\nu_1}^2}{E_1} + \frac{1-{\nu_2}^2}{E_2}.
\end{align}
Elastic modulus and Poisson's ratio for the first material are denoted with $E_1$
and $\nu_1$, and with $E_2$ and $\nu_2$ for the second material. The pressure
distribution between the bodies along the contact surface is semi-elliptical,
i.e.\ of the form
\begin{equation}
p(x) = \begin{cases}
    p_0 \sqrt{1-\frac{x^2}{b^2}}; & |x| \leq b \\
    0; & \text{otherwise}
  \end{cases}, \qquad  p_0 = \sqrt{\frac{PE^*}{\pi R}}.
\end{equation}
A problem can be reduced to two dimensions using plane stress assumption.  A
special case of this problem is when $E_1 = E_2$, $\nu_1 = \nu_2$ and $R_2 \to
\infty$, describing a contact of a cylinder and a half plane. This is the
second numerical example tackled in this paper. The setup is ideal for testing
the refinement, since a pronounced difference in behaviour of numerical
solution near the contact in comparison to the rest of the domain is expected.

A displacement field $\vec{u}$
satisfying~\eqref{eq:navier} on $(-\infty, \infty)
\times (-\infty, 0)$ with boundary conditions
\begin{align}
  \vec{t}(x, 0) = -p(x)\vec{\jmath} \\
  \lim_{x, y\to\infty} \vec{u}(x, y) = 0.
\end{align}
is sought. Vector $\vec{t}$ represents traction force on the surface and
$\vec{\jmath} = (0, 1)$ the upwards direction.  Analytical solution for internal
stresses in the plane in general point $(x, y)$ is calculated using the method
of complex potentials~\cite{mcewen1949stresses} and the stresses are given
in terms of $m$ and $n$, defined as
\begin{align}
  m^2 &= \frac{1}{2} \left(\sqrt{\left(b^2-x^2+y^2\right)^2+4 x^2 y^2}+b^2-x^2+y^2\right), \\
  n^2 &= \frac{1}{2} \left(\sqrt{\left(b^2-x^2+y^2\right)^2+4 x^2 y^2}-(b^2-x^2+y^2)\right),
\end{align}
where $m=\sqrt{m^2}$ in $n=\operatorname{sgn}(x)\sqrt{n^2}$.
The stresses are then expressed as
\begin{align}
  \sigma_{xx} &= -\frac{p_0}{b}\left[m\left(1 + \frac{y^2 + n^2}{m^2 + n^2}\right)+2y\right] \\
  \sigma_{yy} &= -\frac{p_0}{b}m\left(1 - \frac{y^2 + n^2}{m^2 + n^2}\right) \\
  \sigma_{xy} &= \sigma_{yx} = \frac{p_0}{b}n\left(\frac{m^2 - y^2}{m^2 + n^2}\right).
\end{align}

\begin{figure}[h!]
  \centering
  \includegraphics[width=0.6\textwidth]{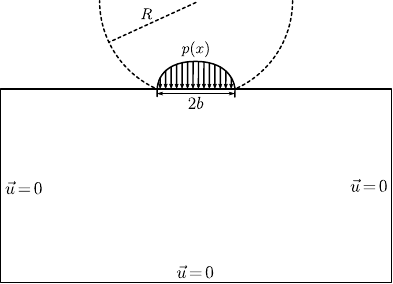}
  \caption{Domain and boundary conditions of considered contact problem.}
  \label{fig:hertz-domain-bc}
\end{figure}

Numerically the problem is solved by truncating the infinite domain to a
rectangle $[-H, H] \times [-H, 0]$ for large enough $H$ and setting the
essential boundary conditions $\vec{u} = 0$ everywhere but on the top boundary.
The top boundary has a traction boundary condition with normal traction given by
$p(x)$ and no tangential traction. An illustration of the problem domain along
with the boundary conditions is given in Figure~\ref{fig:hertz-domain-bc}.

\begin{figure}[h!]
  \centering
  \includegraphics[width=0.7\textwidth]{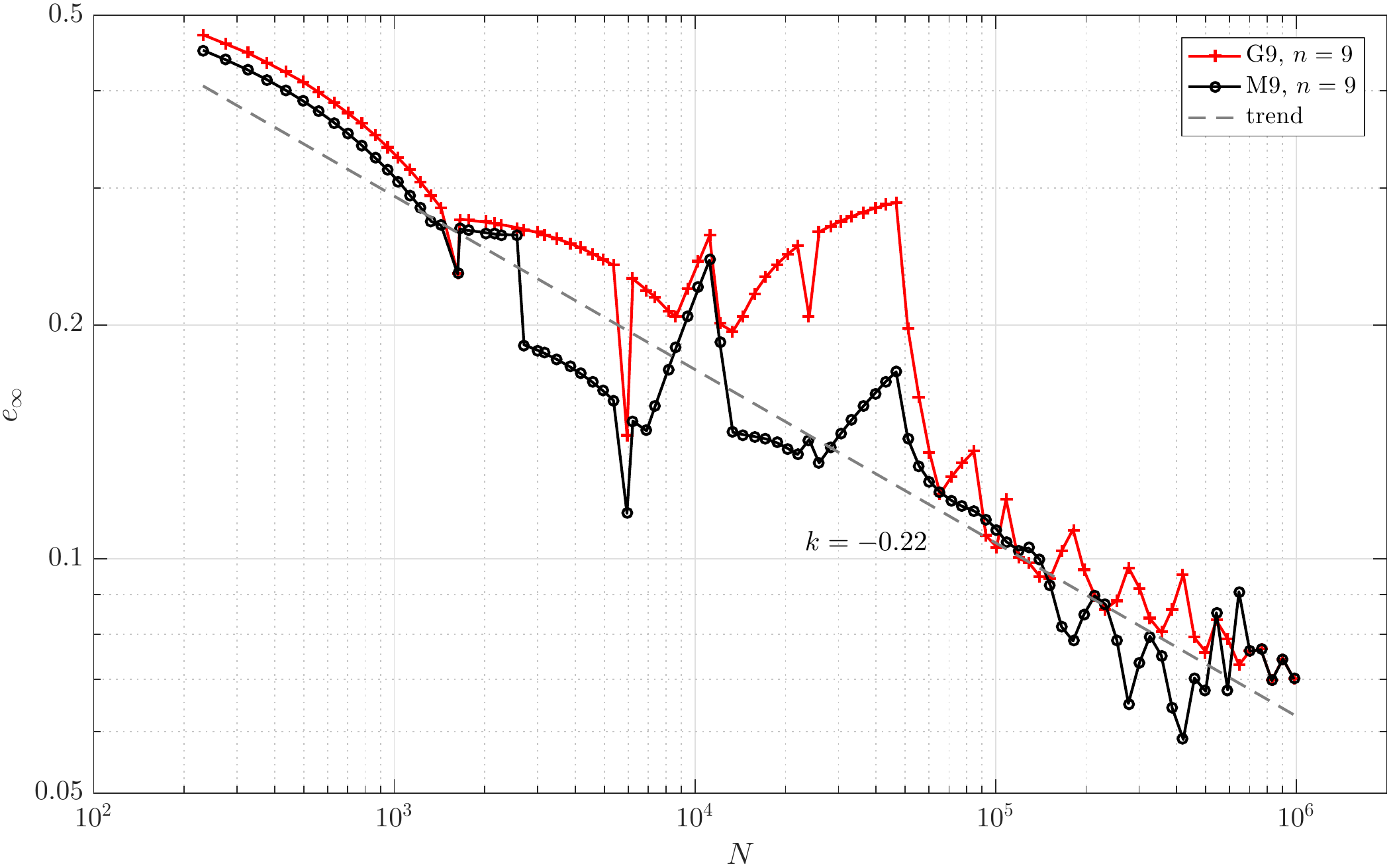}
  \caption{Convergence of MLSM when solving the described Hertzian contact problem.}
  \label{fig:hertz-convergence}
\end{figure}

The described contact problem is solved numerically and the error is measured
between calculated and given stresses in relative $L_\infty$ norm as before,
using
\begin{equation*}
e_\infty = \max_{x\in X} \{\max\{|\sigma_{xx}(x)-\hat{\sigma}_{xx}(x)|,
|\sigma_{yy}(x)-\hat{\sigma}_{yy}(x)|,
|\sigma_{xy}(x)-\hat{\sigma}_{xy}(x)| \}\} / p_0
\end{equation*}
as an error indicator. Values $P = \unit[543]{N/m}$, $E_1 = E_2 =
\unit[72.1]{GPa}$, $\nu_1 = \nu_2 = 0.33$, $R_1 = R = \unit[1]{m}$ were chosen
for the physical parameters of the problem. These values yield contact
half-width $b = \unit[0.13]{mm}$ and peak pressure $p_0 = \unit[2.6]{MPa}$.  A
value of $H = \unit[10]{mm}$ for domain height is chosen, approximately 38 times
greater than width of the contact surface. Convergence of the method
is shown in Figure~\ref{fig:hertz-convergence}.

It is clear that the convergence of the method is very irregular and slow. This
is to be expected, as $N = 10^6$ means only approximately 30 nodes positioned
within the contact surface, and that naturally leads to large changes as a
change of a single node bears a relatively high influence. Another problem is
that the boundary conditions are only continuous, and exhibit no higher
regularity, not even Lipschitz continuity. The accuracy of the approximation
may seem bad, but is in fact comparable to the cantilever beam case. Using the
comparable value of $N = 30\cdot 15 = 450$ nodes in the contact area $[-b, b]
\times [-b, 0]$ it can be seen from Figure~\ref{fig:cantilever-beam-conv} that
the approximation using this number of nodes in cantilever beam case achieved
very similar results.

\begin{figure}[h!]
  \centering
  \includegraphics[width=0.7\textwidth]{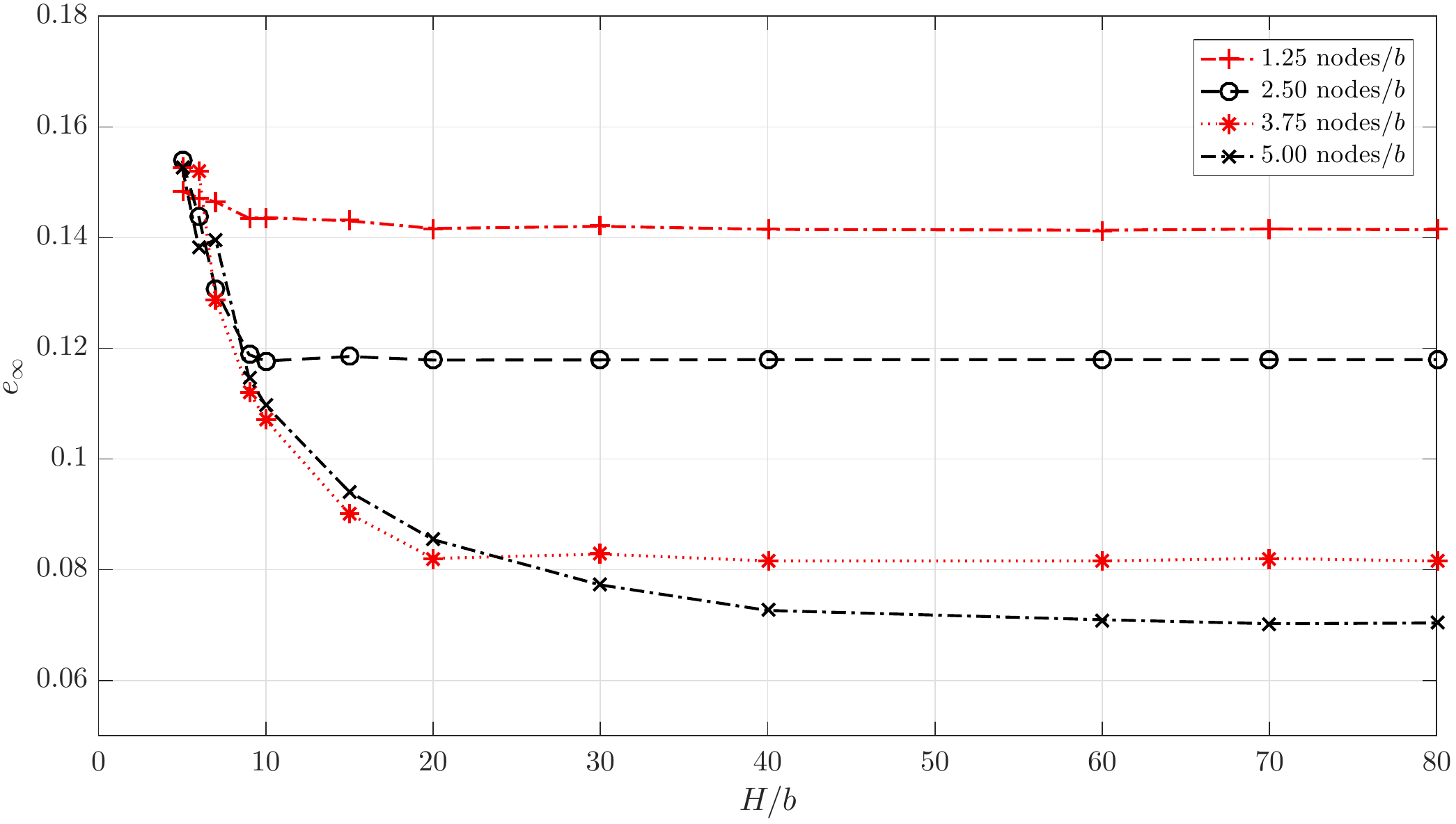}
  \caption{Total error with respect to domain size at different discretization densities.}
  \label{fig:hertz-domain-too-small}
\end{figure}

The total error of the approximation is composed of two main parts, the
truncation error due to the non-exact boundary conditions and the discretization
error, due to solving a discrete problem instead of the continuous one. First, we
analyse the total error in terms of domain height $H$. A graph showing the total
error with respect to domain height $H$ is shown in
Figure~\ref{fig:hertz-domain-too-small}.

The total error decreases as domain height increases, regardless of the
discretization density used. However, as soon as truncation error becomes lower
than discretization error, increasing the height further yields little to no
gain in total error. The higher the discretization density is, the later this
happens. When convergence of a method stops or significantly decreases in
order, an error limit imposed by the truncation error was reached.

It soon becomes impossible to uniformly increase discretization density due to
limited resources, and the immediate solution is to refine the discretization in
the contact area with the h-refinement algorithm introduced in
section~\ref{sec:h-ref}. A domain of height  $H = \unit[1]{m} \approx 75\,000b$
is chosen. Primary refinement is done in rectangle areas of the form
 \[
  [-hb, hb] \times [0, hb], \text{ for } h \in \{1000, 500, 200, 100, 50, 20, 10, 5, 4, 3, 2\},
\]
and secondary refinement around points $\pm b$ on the surface is done in
rectangle areas
\[
[c-hb, c+hb] \times [-hb, 0], \text{ for } c = \pm b \text{ and } h \in \{0.4, 0.3, 0.2, 0.1, 0.05, 0.0025\}.
\]

\begin{figure}[h!]
  \centering
  \begin{subfigure}[t]{0.49\textwidth}
    \includegraphics[width=\textwidth]{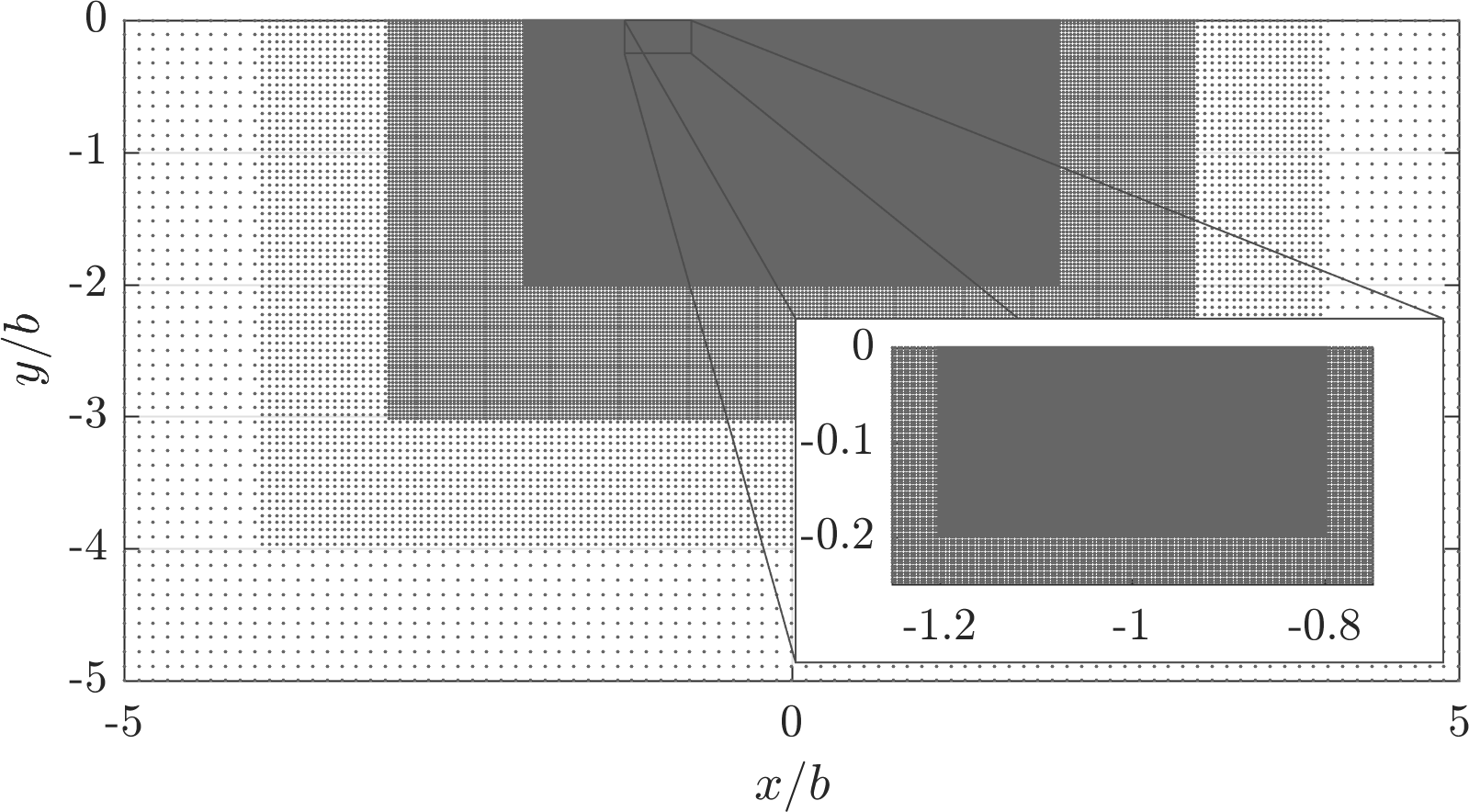}
    \caption{A part of the refined domain.}
    \label{fig:hertz-refined-domain}
  \end{subfigure}
  \begin{subfigure}[t]{0.49\textwidth}
    \includegraphics[width=\textwidth]{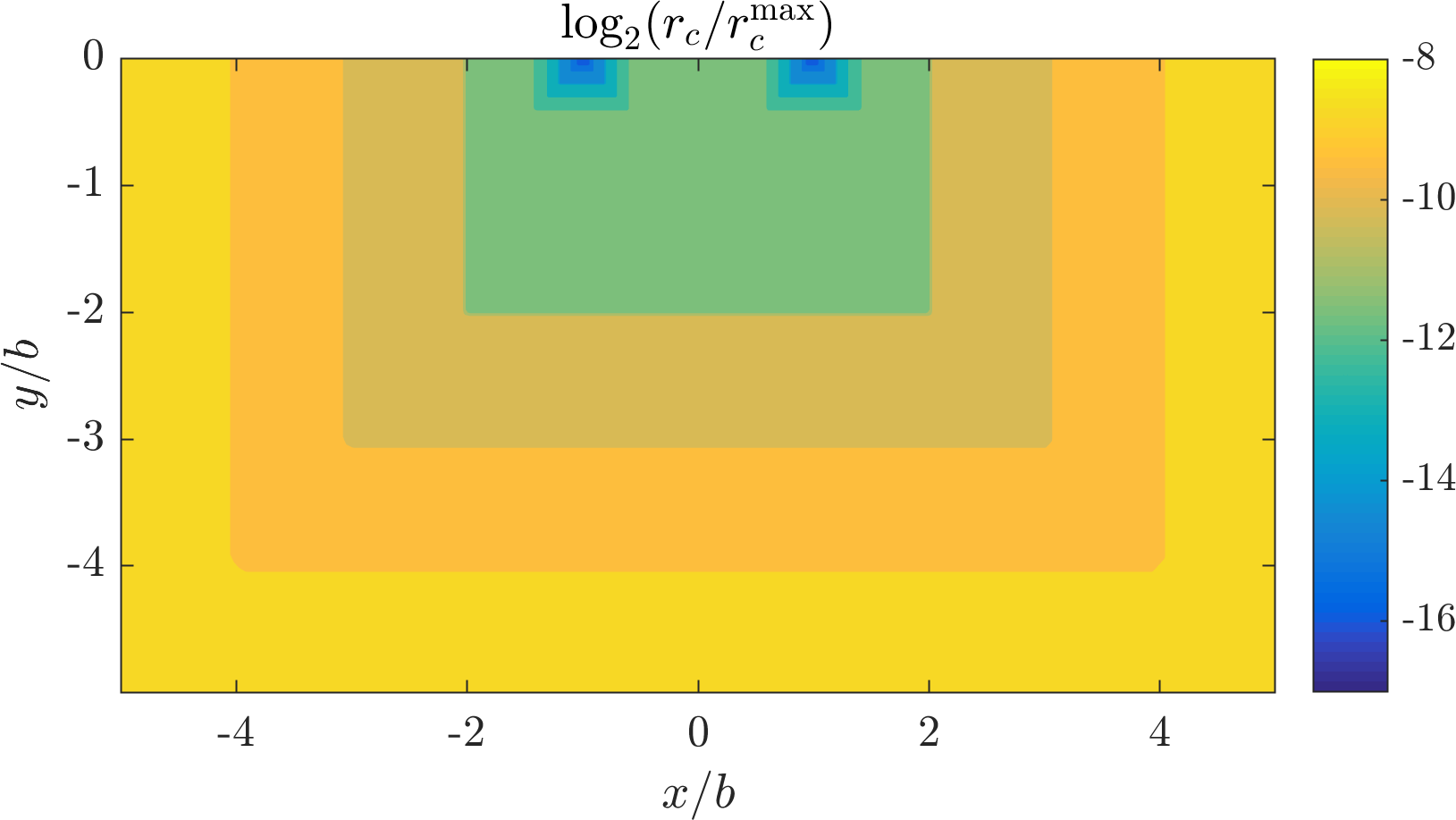}
    \caption{Discretization density of the refined domain.}
    \label{fig:hertz-refined-domain-density}
  \end{subfigure}
  \caption{An example of 17-times refined domain used in solution of
    the described Hertzian contact problem.}
  \label{fig:hertz-refined-domain-together}
\end{figure}

\begin{figure}[h]
  \centering
  \includegraphics[width=0.7\textwidth]{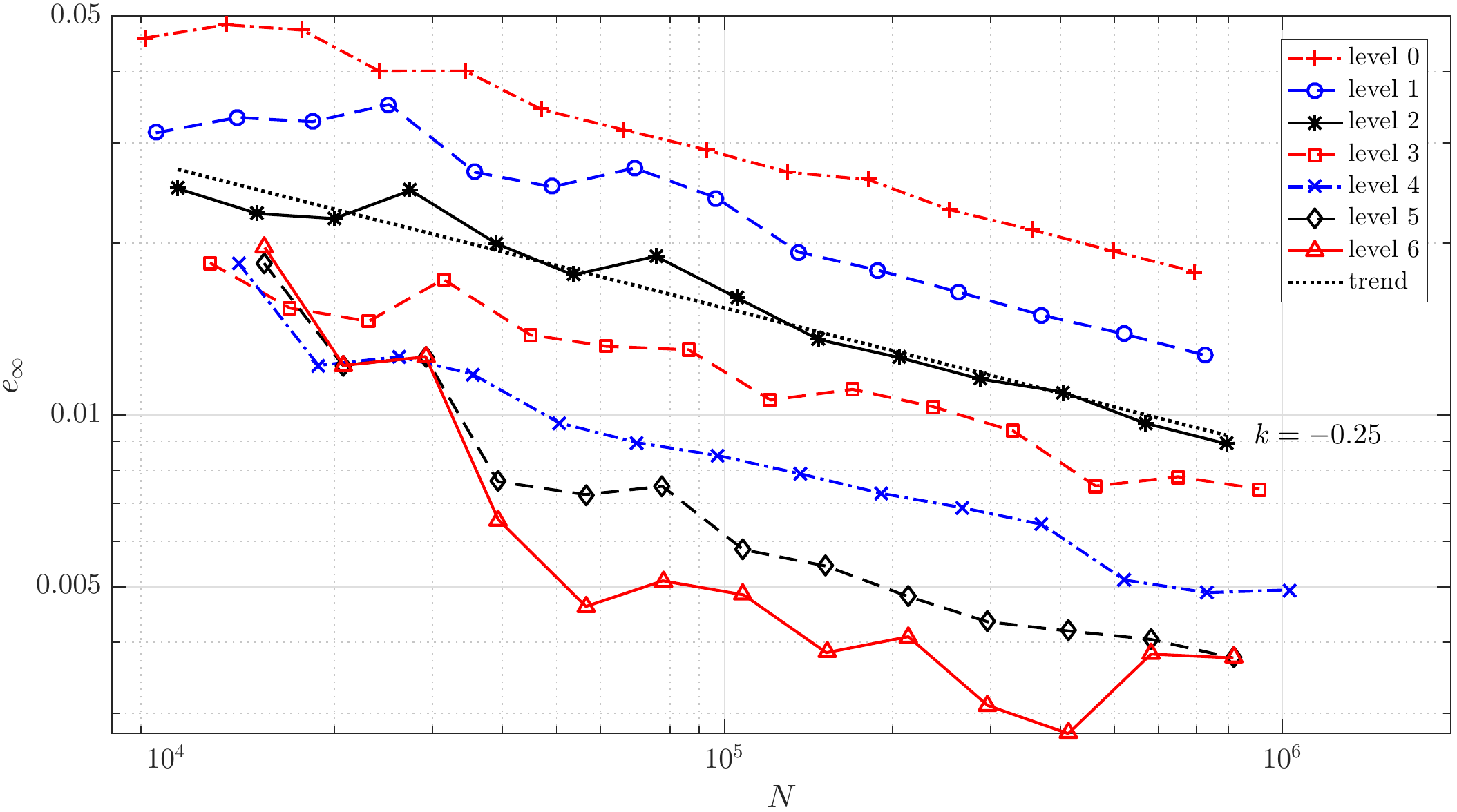}
  \caption{Convergence of MLSM at different refine levels.}
  \label{fig:hertz-refined-convergence}
\end{figure}

The refined domain as described above is shown in
Figure~\ref{fig:hertz-refined-domain-together}.
This domain was used to solve the considered contact problem.
Different levels of secondary refinement were tested to prove that refinement
helps with accuracy. Convergence of the method on the refined domain is shown in
Figure~\ref{fig:hertz-refined-convergence}.

\begin{figure}[h!]
  \centering
  \includegraphics[width=0.70\textwidth]{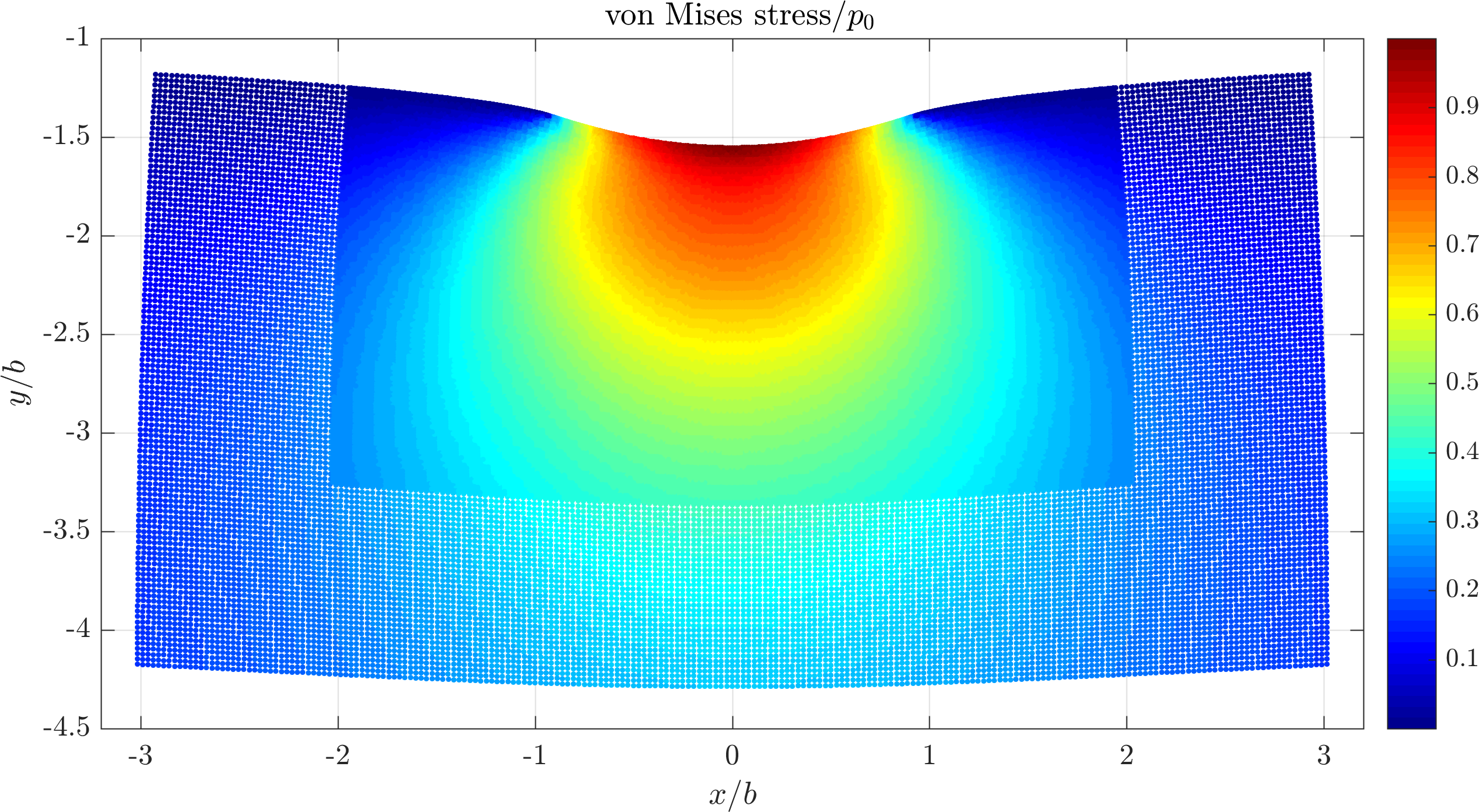}
  \caption{Numerical solution
    of the described Hertzian contact problem. Note that for the
    sake of visibility the displacements are multiplied by factor of $5\cdot 10^3$.}
  \label{fig:hertz-solution}
\end{figure}

Comparing Figure~\ref{fig:hertz-refined-convergence} to
Figure~\ref{fig:hertz-convergence}, it can be seen that refinement greatly
improves the accuracy of the method. Using $N = 10^6$ nodes without refinement
yields worse results than $N = 10^4$ nodes with only primary refinement.  Each
additional level of secondary refinement helps to decrease the error even
further while keeping the same order of convergence. A solution of the problem
on the final mesh is shown in Figure~\ref{fig:hertz-solution}.

\section{Conclusions}
A MLSM solution of a linear elasticity problem on regular and irregular domains
with a refined nodal distribution of two different numerical examples is
presented in this paper. The method is analysed in terms of accuracy by
comparison against available closed form solutions and by comparison against weak form EFG method. The convergence of the
method is evaluated with respect to the number of computational nodes, selection
of different basis functions, different refinement strategies and different
boundary conditions. MLSM is also analysed from complexity point of view, first,
theoretically, and then also experimentally by timing the computer execution
time of all main blocks of the method.
It is clearly demonstrated that the method is accurate and stable. Furthermore, it
is demonstrated that nodal adaptivity is mandatory when solving contact problems
in order to obtain accurate results and that the proposed MLSM method can handle
extensive refinement with the smallest internodal distance being $2^{17}$ times
smaller than the initial one. It is also demonstrated that proposed MLSM
configuration can handle computations in complex domains.

In our opinion the presented meshless setup can be used, not only to solve
academic cases with the sole goal to show excellent convergences, but also in
more complex engineering problems. The C++ implementation of presented MLSM is
freely available at~\cite{utils_web}.

In future work we will continue to develop a meshless solution of a contact
problem with a final goal to simulate a crack propagation due to the fretting
fatigue~\cite{pereira2016convergence} in a general 3D domain with added
p-adaptivity to treat singularities near the crack tip.

\section*{Acknowledgement}
The authors would like to acknowledge the financial
support of the Research Foundation Flanders (FWO), The Luxembourg National
Research Fund (FNR) and Slovenian Research Agency (ARRS) in the framework of the
FWO Lead Agency project: G018916N Multi-analysis of fretting fatigue using
physical and virtual experiments and the ARRS research core funding No.
P2-0095.

\bibliographystyle{unsrt}
\bibliography{mybibfile}

\end{document}